\documentclass[reqno]{amsart}

\usepackage{graphicx}
\usepackage{amsmath,amssymb}
\usepackage{amsthm}
\usepackage{xcolor}
\usepackage{mathrsfs}
\usepackage{paralist}
\usepackage{multirow}
\usepackage[T1]{fontenc}
\usepackage[utf8]{inputenc}
\usepackage[english]{babel}
\DeclareMathAlphabet\mathbfcal{OMS}{cmsy}{b}{n}
\usepackage[figurename=Figure]{caption}
\usepackage{lineno}
\usepackage{bbm}
\usepackage[normalem]{ulem}
\usepackage[percent]{overpic}
\usepackage{color}
\usepackage{tikz}
\usepackage{yhmath}

\usepackage{wrapfig}
\usepackage{float}
\usepackage{booktabs}
\usepackage{makecell}
\usepackage{subcaption}
\usepackage{comment}

\usepackage{hyperref}
\hypersetup{
    colorlinks=true,     
    linkcolor=red,       
    citecolor=green,     
    filecolor=magenta,   
    urlcolor=blue        
}

\newcommand{\del}[1]{}

\usepackage{authblk}

\usepackage{amsaddr}
\usepackage{siunitx}

\usepackage[a4paper,bindingoffset=0.0in,%
            left=1.0in,right=1.0in,top=1.in,bottom=1.in,%
            footskip=.4in]{geometry}

\theoremstyle{remark}
\newtheorem{remark}{Remark}[section]
\title[Efficient and Accurate Surrogate Modeling via Aggregation and ROMs]{Efficient and Accurate Surrogate Modeling of Turbulent Flows via Space-Dependent Aggregation and Reduced Order Models}

\author{Piero Zappi$^{1}$, Anna Ivagnes$^{1}$, Niccolò Tonicello$^{1}$ and Gianluigi Rozza$^{1}$}
\address{$^1$ mathLab, Mathematics Area, SISSA, via Bonomea 265, I-34136 Trieste, Italy}

\begin{document}
\begin{abstract}
Reynolds-Averaged Navier–Stokes (RANS) models are widely used for turbulent flow simulations due to their computational efficiency, but their accuracy strongly depends on the selected turbulence closure and may vary across the flow domain. Space-dependent model aggregation has been shown to improve RANS predictions by combining multiple turbulence models, although at the cost of repeated high-fidelity simulations.
The first novelty of this work is a unified framework that combines different turbulence models, space-dependent aggregation, and non-intrusive reduced order models to achieve both accuracy and efficiency. Two aggregation pipelines are proposed: a Mixed FOM–ROM (MFR) approach, where a reduced order model is trained on aggregated RANS solutions, and a Mixed-ROM (MR) approach, which directly aggregates multiple reduced order models built on top of different RANS full-order models. The second novelty is that the aggregation weights are learned via a neural-network that provides smooth, space-continuous weights and improves generalization with respect to standard weighting techniques.
The resulting surrogate models are validated on the two-dimensional periodic hill benchmark and on the flow over a height-dependent bump, demonstrating improved accuracy over individual RANS and ROM predictions at near real-time computational cost.
\end{abstract}

\maketitle

\section{Introduction}
Computational Fluid Dynamics (CFD) plays a central role in a wide range of scientific and engineering disciplines, including aerospace \cite{slotnick2014cfd} and automotive engineering \cite{aultman2022evaluation}, energy systems \cite{tominaga2024cfd}, environmental modeling, and geophysical flows \cite{blocken2007cfd}. 
By enabling the numerical solution of the governing equations of fluid motion, CFD has become an indispensable tool for the analysis, design, and optimization of complex systems. Continuous advances in numerical methods and high-performance computing have significantly expanded the range of problems that can be tackled. Nevertheless, the simulation of turbulent flows remains one of the most demanding challenges in CFD, due to the wide range of interacting spatial and temporal scales involved.

From the standpoint of turbulence modeling, three main levels of fidelity are commonly adopted: Direct Numerical Simulation (DNS), Large Eddy Simulation (LES), and Reynolds-Averaged Navier–Stokes (RANS) approaches. DNS resolves all turbulent scales and provides the highest level of accuracy, but its computational cost grows prohibitively with the Reynolds number, making it unfeasible for most practical applications. LES offers a compromise by resolving the large energy-containing structures while modeling the smaller scales, but it remains too expensive for many industrial and multi-query scenarios \cite{moin1997tackling}. As a consequence, RANS models are still the workhorse of industrial CFD, due to their relatively low computational cost and robustness.

However, RANS models often suffer from well-known limitations in predictive accuracy for complex turbulent flows involving, for example, separation, strong adverse pressure gradients, or geometry-induced effects \cite{cinnella2024data}. Different turbulence closures may lead to substantially different flow predictions \cite{xiao2019quantification}, and models are often designed and calibrated for specific classes of flows. As a result, no single RANS turbulence model can be expected to perform optimally across all flow regimes or even across all regions of the same computational domain. This has motivated the development of model aggregation strategies, whose objective is to combine the predictions of multiple models in order to leverage their complementary strengths and mitigate individual weaknesses.

Model aggregation has a long history in statistics and forecasting, where techniques such as the Exponentially Weighted Average (EWA) were introduced to combine multiple predictors into a single, more robust estimate \cite{MA_1_EWA,MA_2,MA_3}. In the CFD context, these aggregation strategies have been adapted to turbulence modeling, initially through spatially uniform convex combinations of RANS solutions. While effective, these approaches neglect the inherently local nature of turbulence modeling errors. More recent developments have therefore focused on space-dependent aggregation, inspired by Mixture-of-Experts methods \cite{Mixtures_of_experts, yuksel2012twenty}, in which the contribution of each model varies across the domain according to local flow features. The space-dependent aggregation framework combines the robustness of model averaging with the adaptivity of local expert selection and has been shown to significantly improve RANS predictions \cite{aggregation_Cinnella1,aggregation_Cinnella2,oulghelou2025machine}. However, a key limitation of these methods is that each aggregated prediction typically requires running several RANS simulations, which can still be computationally expensive in multi-query settings.

This observation naturally motivates the development of a further layer of computational efficiency via surrogate modelling. One possible option is represented by Reduced Order Models (ROMs) for RANS, whose goal is to further accelerate simulations while retaining acceptable accuracy. Data-driven and non-intrusive ROMs, in particular, have proven effective in exploiting the intrinsic low-dimensional structure of RANS solution manifolds~\cite{maulik2021reduced,saetta2022machine,racca2023predicting,macedo2024using,aggregation_anna,tonioni2025data,fung2025vector,halder2026reduced}.

The main idea of this work is to combine space-dependent model aggregation, different turbulence models, and reduced order methods, to simultaneously improve accuracy and drastically reduce computational costs. We introduce two complementary aggregation pipelines. The first, referred to as Mixed FOM–ROM (MFR), constructs a reduced order model downstream of aggregated full-order RANS solutions. The second, denoted as Mixed-ROM (MR), instead aggregates multiple ROMs, each associated with a different RANS turbulence model, directly at the reduced-order level. In both cases, the final surrogate models operate entirely in an online-efficient setting, eliminating the need for new RANS simulations at prediction time.

A central component of both pipelines is the computation of space-dependent aggregation weights. A common strategy in literature is  space-dependent EWA aggregation with a K-Nearest Neighbors (KNN) regression to predict the weights for unseen configurations, resulting in the MFR-KNN and MR-KNN models.
In this work, we introduce the MFR-ANN and MR-ANN models, which employ a feed-forward Artificial Neural Network (ANN) that directly outputs spatially continuous weights. This enables improved generalization and accuracy for unseen configurations with respect to standard approaches. The ANN-based approach removes the need for kernel-based weight reconstruction and allows for smoother, fully continuous weight fields.

The proposed methodologies are assessed on the test cases of the two-dimensional incompressible flow over periodic hills, and of the flow over a bump. Both test cases involve geometric parametrization and are widely regarded as benchmark configurations for turbulence modeling because of their intricate flow separation and reattachment dynamics~\cite{wu2023enhancing, Xiao, Breuer, Frohlich, Rapp, CinnellaLES, periodic_hill_extra1, periodic_hill_extra2, periodic_hill_extra3, bump1_LES, bump2, bump3, bump4, bump5, bump6}. 

The remainder of the paper is organized as follows. Section \ref{sec:methods} introduces the methodological background, including RANS turbulence modeling, reduced order modeling, and space-dependent aggregation techniques. Section \ref{sec:aggregation-novelty} presents the proposed aggregation pipelines and the two weighting strategies we investigated. Section \ref{sec:results} reports the numerical results and a detailed comparison of accuracy and computational efficiency. Finally, Section \ref{sec:conclusions} concludes the paper and outlines possible directions for future work.

\section{Methodological background}
\label{sec:methods}
In the general framework we propose, three main methods play a central role in the subsequent sections: the RANS approach (Subsection \ref{sec:RANS}), non-intrusive Reduced Order Models (Subsection \ref{sec:ROMs}), and model aggregation (Subsection \ref{sec:aggregation}). The RANS equations provide the physical foundation for the simulation of turbulent flows, while Reduced Order Models and aggregation act in a complementary manner to enhance computational efficiency and predictive accuracy, respectively. In the following sections, we not only introduce each of these topics, but also discuss how their interplay leads to a self-contained framework for fast and accurate prediction of turbulent flow statistics.

\subsection{RANS approach}\label{sec:RANS}
\par The RANS formulation is obtained by applying an averaging operator to the Navier-Stokes equations. The cornerstone of this approach is the \textit{Reynolds decomposition} \cite{Reynolds}, which assumes that any instantaneous flow quantity can be written as the sum of a mean component and a zero-mean fluctuating component. For a generic space- and time-dependent field $\phi(\mathbf{x},t)$, this decomposition reads:
\begin{equation*}    \phi(\mathbf{x},t)=\overline{\phi}(\mathbf{x},t)+\phi'(\mathbf{x},t),
\end{equation*}
\par where $\overline{\phi}(\mathbf{x},t)$ denotes the mean field and $\phi'(\mathbf{x},t)$ the fluctuating component. In the most general setting, the mean operator is defined as an ensemble average. In many practical applications, however, the presence of symmetries or an expected statistically steady flow allows this definition to be simplified to time averaging and/or spatial averaging \cite{wilcox1998turbulence}.

\par To characterize different flow regimes, we consider the dimensionless Reynolds number $\mathrm{Re}$, namely the ratio between inertial and viscous effects:
\begin{equation*}
    \mathrm{Re}=\displaystyle{\frac{UL}{\nu}},
\end{equation*}
\par where $U$ and $L$ are the characteristic case-dependent velocity and length scales, respectively, and $\nu$ is the kinematic viscosity. For a given problem, it is generally possible to identify a critical value of $\mathrm{Re}$ above which inertial effects prevail over viscous dissipation, causing the flow to transition from a laminar to a turbulent regime.
\par Let $\mathbf{u}(\mathbf{x},t)$ denote the velocity field and $p(\mathbf{x},t)$ the normalized pressure. Substituting the Reynolds decomposition into the Navier--Stokes equations and taking the ensemble average leads to the RANS equations (written in the Einstein notation):
\begin{equation}
    \begin{cases}
    \displaystyle{
    \frac{\partial \overline{u}_i}{\partial x_i}=0},\\[6pt]
    \displaystyle{
    \frac{\partial \overline{u}_i}{\partial t}
    +\overline{u}_j\frac{\partial \overline{u}_i}{\partial x_j}
    =-\frac{\partial \overline{p}}{\partial x_i}
    +\frac{1}{\mathrm{Re}} \frac{\partial^{2} \overline{u}_{i}}{\partial x_{j} \partial x_{j}}
    - \frac{\partial \mathcal{R}_{ij}}{\partial x_{j}} ,
    }
    \end{cases}
    \label{RANSeq}
\end{equation}
where $\mathcal{R}_{ij}=\overline{u'_i u'_j}$ is the Reynolds stress tensor (RST). Since the velocity fluctuations are not available in RANS simulations, closure models must be introduced to approximate the Reynolds stresses, which constitutes the central challenge of turbulence modeling.

\par As a covariance tensor, $\mathcal{R}_{ij}$ is symmetric and positive semi-definite, a property commonly referred to as the \emph{realizability requirement}. The vast majority of turbulence models target the modeling of the full Reynolds stress tensor although a modelisation for its divergence is also possible. The former approach, however, is generally preferred as it guarantees conservation of momentum over any arbitrary control volume.

A broad range of turbulence models has been developed over the years. Despite sustained research efforts over several decades, no single modeling approach has proven to be universally reliable or consistently adequate across all flow regimes \cite{spalart2000strategies}. In particular, RANS models are known to exhibit significant limitations when applied to complex flow configurations characterized by turbulence non-equilibrium, strong gradients, flow separation, shocks, and pronounced three-dimensional effects. Moreover, multiple sources of uncertainty influence the reliability of RANS predictions, which further hinders the design of turbulence models that are robust and transferable across different turbulent flow regimes \cite{duraisamy2019turbulence,xiao2019quantification}.
As a result, flow-specific tuning and ad hoc corrections remain an essential component of many RANS simulations \cite{spalart2015philosophies}. The ongoing challenge in turbulence modeling is therefore to preserve the low computational cost and robustness of RANS approaches while improving their predictive accuracy toward that of high-fidelity simulations.

Reynolds-stress-based turbulence models require prescribing a constitutive relation for $\mathcal{R}_{ij}$ in terms of the mean flow quantities. The most widely used class of models, commonly referred to as \emph{linear eddy viscosity models} (LEVMs), is based on the Boussinesq hypothesis (see, e.g., \cite{wilcox1998turbulence}). This assumption postulates that the anisotropic part of the Reynolds stress tensor behaves analogously to the viscous stress tensor of a Newtonian fluid and is therefore proportional to the local mean rate-of-strain tensor $\overline{S}_{ij}$:
\begin{equation}
    -\mathcal{R}_{ij}
    = 2\nu_t \overline{S}_{ij}
    - \frac{2}{3} k \delta_{ij},
    \quad
    \text{with} \quad
    \overline{S}_{ij}
    = \frac{1}{2}
    \left(
    \frac{\partial \overline{u}_{i}}{\partial x_{j}}
    + \frac{\partial \overline{u}_{j}}{\partial x_{i}}
    \right),
    \quad
    k = \frac{1}{2}\overline{u'_i u'_i}.
\end{equation}
\par Here, $\nu_t$ denotes the turbulent (or \emph{eddy}) viscosity and $k$ is the turbulent kinetic energy. The eddy viscosity accounts for the influence of turbulence on the mean flow and represents a non-physical, turbulence-induced analogue of the molecular viscosity. Unlike the latter, $\nu_t$ is not a material property of the fluid but depends on the local flow conditions. Different turbulence models provide different transport equations or algebraic relations for computing $\nu_t$, thereby leading to different closures of the RANS equations. In this work, we focus on the following widely used eddy viscosity models: the \emph{Spalart-Allmaras} (SA) model \cite{SA_article}, the $k\varepsilon$ model \cite{kappaepsilon}, the $k\omega$ model \cite{kappaomega}, and the $k\omega$ SST model \cite{kappaomegasst}.

\subsection{Non-intrusive Reduced Order Models}\label{sec:ROMs}
\par Reduced Order Models are a broad class of surrogate models designed to efficiently approximate high-dimensional systems by projecting their dynamics onto low-dimensional subspaces, thereby enabling significantly faster simulations \cite{bookBenner,aromabook,bookROMS1,bookROMS2,bookROMS3,bookROMS4,bookROMS5}. Their efficiency relies on the so-called \textit{offline-online} procedure. The offline phase, performed only once, includes all computationally demanding tasks, such as setting up the full-order model (FOM) and running a sufficiently large set of high-fidelity simulations over a range of parameter values to build the training dataset and ensure that the solution manifold adequately represents the test case. The online phase, in contrast, is designed to be inexpensive and may be executed repeatedly; in this stage, the ROM provides reliable approximations for previously unseen parameter configurations without requiring any additional costly simulations.
\par In this work, we concentrate on \textit{data-driven} \textit{non-intrusive} ROMs, which depend exclusively on data and operate independently of the underlying physics. Unlike intrusive approaches, they do not require any knowledge of either the governing equations or the numerical schemes employed to simulate the physical system. This makes non-intrusive ROMs particularly advantageous in scenarios where the governing equations are unavailable but high-fidelity data exist (for example, from experiments), or in cases where access to the full-order models is restricted, such as when proprietary commercial solvers are employed.
\par The full-order solutions, referred to as \textit{snapshots}, represent the fields of interest corresponding to different parameter configurations. For a generic field $s$, the snapshot associated with parameters $\boldsymbol{\mu}_i\in\mathbb{R}^{p}$ is denoted by $\mathbf{s}_i=\mathbf{s}(\boldsymbol{\mu}_i)$. These snapshots are collected into the matrix $\mathbf{S}$:
\begin{equation}
\mathbf{S} = \left[
\begin{array}{cccc}
| & | & & | \\
\mathbf{s}_1(\mathbf{x}) & \mathbf{s}_2(\mathbf{x}) & \cdots & \mathbf{s}_N(\mathbf{x})\\
| & | & & |
\end{array}
\right]
\in \mathbb{R}^{n_{\mathrm{dof}} \times N_{\mathrm{train}}},\label{SnapMatrix}
\end{equation}
\par where $n_{\mathrm{dof}}$ is the number of degrees of freedom, which remains fixed across all full-order simulations to ensure a consistent representation. $N_{\mathrm{train}}$ is the number of training snapshots.
\par The offline phase of a non-intrusive ROM consists of two fundamental stages. First, in the \textit{reduction} step, a suitable low-dimensional representation of the high-fidelity data is constructed by performing a compression of the snapshots onto a reduced space. Second, the \textit{approximation} step involves learning a mapping between the problem's parameters and the reduced variables, typically through interpolation or regression approaches. During the online stage, the learned mapping is queried for a new parameter value $\boldsymbol{\mu}^*$ to predict its reduced representation, which is then lifted back to the original space to obtain an approximation $\widetilde{\mathbf{s}}(\boldsymbol{\mu}^*)$. This procedure is extremely efficient and enables fast and accurate predictions for previously unseen configurations without requiring access to the high-fidelity model, and thereby providing substantial computational savings.

\par The reduction step consists of applying a linear or nonlinear mapping $\mathscr{R}$ to compress the snapshots' matrix $\mathbf{S}$ (\ref{SnapMatrix}) into a reduced space of dimension $r$. Nonlinear reduction is particularly suitable for advection-dominated problems, as the nonlinear mapping allows for the representation of complex flow features that linear approaches often fail to capture effectively. In this work, we employ autoencoders (AEs) as nonlinear reduction tools. Autoencoders are a class of neural networks that have gained popularity as nonlinear dimensionality reduction techniques \cite{aggregation_anna,AE1,AE2,AE3,AE4,AE5,AE6,AE7,AE8} due to their distinctive architecture, which inherently performs both the encoding of high-dimensional data into a latent space and the subsequent reconstruction back to the original dimensionality. An AE is composed of two components: an \textit{encoder} $\mathcal{\epsilon}:\mathbb{R}^{n_{\mathrm{dof}}}\rightarrow\mathbb{R}^{r}$, which maps the high-dimensional input data into a reduced space, and a \textit{decoder} $\mathcal{\delta}:\mathbb{R}^{r}\rightarrow\mathbb{R}^{n_{\mathrm{dof}}}$, which reconstructs the full-dimensional representation from the latent encoding. Autoencoders can be constructed using either convolutional \cite{lee2020model} or fully connected layers \cite{wang2014generalized, wang2016auto}, and numerous variants have been proposed in the literature of model order reduction of parametric PDEs. Their applications extend beyond model order reduction to areas such as anomaly detection \cite{chen2018autoencoder}, image processing \cite{balle2016end}, and neural machine translation \cite{wu2014proceedings}. In the present study, both the encoder and the decoder are implemented as dense feed-forward neural networks. Within the framework of model order reduction technique, AEs have shown particularly significant success in treating complex fluid dynamics problems involving turbulent flows \cite{maulik2021reduced,saetta2022machine,racca2023predicting,macedo2024using,tonioni2025data,fung2025vector,halder2026reduced}.

\par The approximation step aims to learn a mapping from the parameter space to the latent space, enabling predictions for previously unseen parameter values during the online phase. In this work, Radial Basis Function (RBF) Interpolation \cite{rbf} is employed due to its simplicity, flexibility, and satisfactory predictive accuracy.


\subsection{Model aggregation}\label{sec:aggregation}
\par In this study, we build the model aggregation strategy on the XMA algorithm \cite{aggregation_Cinnella1, aggregation_Cinnella2}, which constructs mixtures of RANS models with the objective of achieving predictive performance superior to that of any individual component. The core idea of such a framework is to combine the outputs of multiple turbulence models through a space-dependent weighted mixture, leveraging their complementary strengths to improve overall accuracy.
\par Let $\mathcal{M}=\{M_1,M_2,\dots,M_{n_M}\}$ denote a discrete set of $n_M$ full-order models, each corresponding to a CFD simulation of the same physical problem but differing in the turbulence closure. The prediction produced by model $M_i$ is denoted by:
\begin{equation*}
    \widetilde{\mathbf{s}}^{(i)}=\widetilde{\mathbf{s}}(\boldsymbol{\eta},M_i),\quad i=1,\dots,n_M,
\end{equation*}
\par where $\boldsymbol{\eta}$ represents the independent input features driving the prediction. The aggregated prediction is defined as a convex combination of the individual model outputs:
\begin{equation}
    \widetilde{\mathbf{s}}^{(mix)}=\sum_{i=1}^{n_M}\omega_i(\boldsymbol{\eta})\widetilde{\mathbf{s}}^{(i)},\label{aggr_pred}
\end{equation}
\par where $\omega_i(\boldsymbol{\eta})$ denotes the weights’ distribution assigned to model $M_i$. These weights are constrained to satisfy the conditions:
\begin{equation}
\omega_i(\boldsymbol{\eta})\in[0,1],\quad\sum_{i=1}^{n_M}\omega_i(\boldsymbol{\eta})=1,\quad\forall\boldsymbol{\eta}.\label{norm_conditions}
\end{equation}
\par This formulation assigns larger weights to models exhibiting higher predictive accuracy, while giving smaller weights to less accurate ones, resulting in a mixture that outperforms its components. Importantly, the aggregation procedure is space-dependent, and the weights vary in space, enabling each individual model to contribute primarily in the regions of the computational domain where it is most accurate.

\section{Aggregation pipelines and weighting strategies}
\label{sec:aggregation-novelty}
\par In this Section, we propose two different aggregation pipelines for combining reduced-order models and aggregation techniques for the simulation of turbulent flows. The objective is to enhance the final predictive accuracy by leveraging the complementary strengths of different turbulence closures through localized aggregation while significantly reduce the computational costs via ROMs as efficient surrogate models.

\subsection{Weighting techniques}
\par The predictive performance of the model mixtures depends on how the outputs of the component models are combined through spatially varying weights, which quantify the local reliability of each model and determine its contribution to the aggregated solution. In this work, two alternative strategies for computing optimal space-dependent weights are investigated.

\par The first technique is a modified version of the space-dependent Model Aggregation (XMA) algorithm \cite{aggregation_Cinnella1, aggregation_Cinnella2}, later adopted in \cite{aggregation_anna}, and is based on the Exponentially Weighted Average (EWA) predictor \cite{MA_1_EWA}. The mixture weights are defined as
\begin{equation}
    \omega_i(\boldsymbol{\eta})=\omega_i(\boldsymbol{\mu})=\frac{g_i(\boldsymbol{\mu})}{\sum_{j=1}^{n_M}g_j(\boldsymbol{\mu})},\quad i=1,\dots,n_M,
\end{equation}
\par where $\boldsymbol{\mu}$ denotes the physical parameters of the simulation and the cost functions $g_i$ are modeled as Gaussians:
\begin{equation}
    g_i(\boldsymbol{\mu})=\exp\!\left(-\frac{1}{2}\frac{(\widetilde{\mathbf{s}}^{(i)}-\mathbf{s})^2}{\varsigma^2}\right).
\end{equation}

\par Here, $\mathbf{s}$ is a reference solution, taken in this work from DNS snapshots. The cost function attains its maximum value when the model prediction matches the reference and decays for increasing discrepancies. The hyperparameter $\varsigma$ controls the sensitivity of the EWA to prediction errors: large values yield nearly uniform weights, while small values promote highly selective weighting. The parameter $\varsigma$ has been obtained as in \cite{MA_1_EWA,aggregation_Cinnella1}. 


\par To generalize the weighting strategy to unseen flow configurations, supervised regressors are trained to learn the mapping between input parameters and the optimal weights.
We consider a number of training parameters for the aggregation $N_{\mathrm{train}}^{(mix)}$.
Given the training pairs $(\boldsymbol{\mu}_{\mathrm{train}}^{(mix)},\omega_j(\boldsymbol{\mu}_{\mathrm{train}}^{(mix)}))$, the weights for new configurations $\boldsymbol{\mu}_{\mathrm{test}}$ are predicted. In this work, $n_M$ independent KNN regressors   \cite{scikit-learn} are used to approximate the cost functions $g_i$, which are subsequently normalized.

In addition to the classical approach based on KNN, in this work a novel strategy is introduced to calculate mixture weights varying spatially that take advantage of the capabilities of artificial neural networks \cite{universalapprox1,universalapprox2} to enhance the accuracy of the aggregated predictions. In this approach, a feed-forward fully connected Artificial Neural Network (ANN) is employed. The ANN takes as input both the spatial coordinates $\boldsymbol{x}$ and the physical parameters $\boldsymbol{\mu}$ of the simulations and outputs directly the normalized weights $\omega_j$, $j=1,\dots,n_M$, associated with each $j$-th component model. A \textit{softmax} layer is added at the end of the network to ensure that such weights satisfy the normalization and positivity constraints. The ANN is trained by minimizing the mean squared error between the aggregated predictions $\widetilde{\mathbf{s}}^{(mix)}$ and the corresponding reference solutions $\mathbf{s}$, according to the loss function:
\begin{equation}
    \mathcal{L}=\displaystyle{\frac{1}{N_{\mathrm{train}}^{(mix)}}}\sum_{i=1}^{N_{\mathrm{train}}^{(mix)}}\bigg((\mathbf{s})_i-(\widetilde{\mathbf{s}}^{(mix)})_i\bigg)^2.
\end{equation}
\par Details regarding the ANN architecture and the associated training hyperparameters are provided in the supplementary material (Section \ref{sec:sup}). The network parameters are initialized using the \textit{He} initialization scheme \cite{He_init} to promote stable convergence and the Adam optimizer \cite{Adam} is used for training.
\par This ANN-based weighting strategy offers several advantages over the Gaussian-based approach.
\begin{enumerate}
    \item It is expected to deliver superior performance by exploiting a richer set of input features, $\boldsymbol{\eta}=\boldsymbol{\eta}(\boldsymbol{x},\boldsymbol{\mu})$, which incorporates both spatial coordinates and physical parameters of the simulations.
    \item It removes the dependence on the hyperparameter $\varsigma$, thereby eliminating the need for a dedicated optimization procedure and simplifying the overall workflow.
    \item The ANN defines a continuous mapping from the input space to the mixture weights, enabling straightforward generalization to previously unseen spatial locations and flow configurations. This property is particularly advantageous in scenarios involving datasets of varying sizes, as commonly encountered in multi-fidelity modeling frameworks.
\end{enumerate}

\begin{remark}
    It is important to remark that the training set used for the ROM may be different from the training set used for the aggregation strategy. Consequently, in general $N_{\mathrm{train}}^{(mix)} \neq N_{\mathrm{train}}$. In particular, in the first test case, we will consider exactly the same training set for ROM and aggregation, while for the second test case, as we have much less data available, $N_{\mathrm{train}}^{(mix)} \neq N_{\mathrm{train}}$.
\end{remark}

\subsection{ROM-aggregation pipelines}
\par
\begin{figure}[!htpb]
     \centering
     \includegraphics[trim={350 0 350 0}, clip, width=0.65\textwidth]{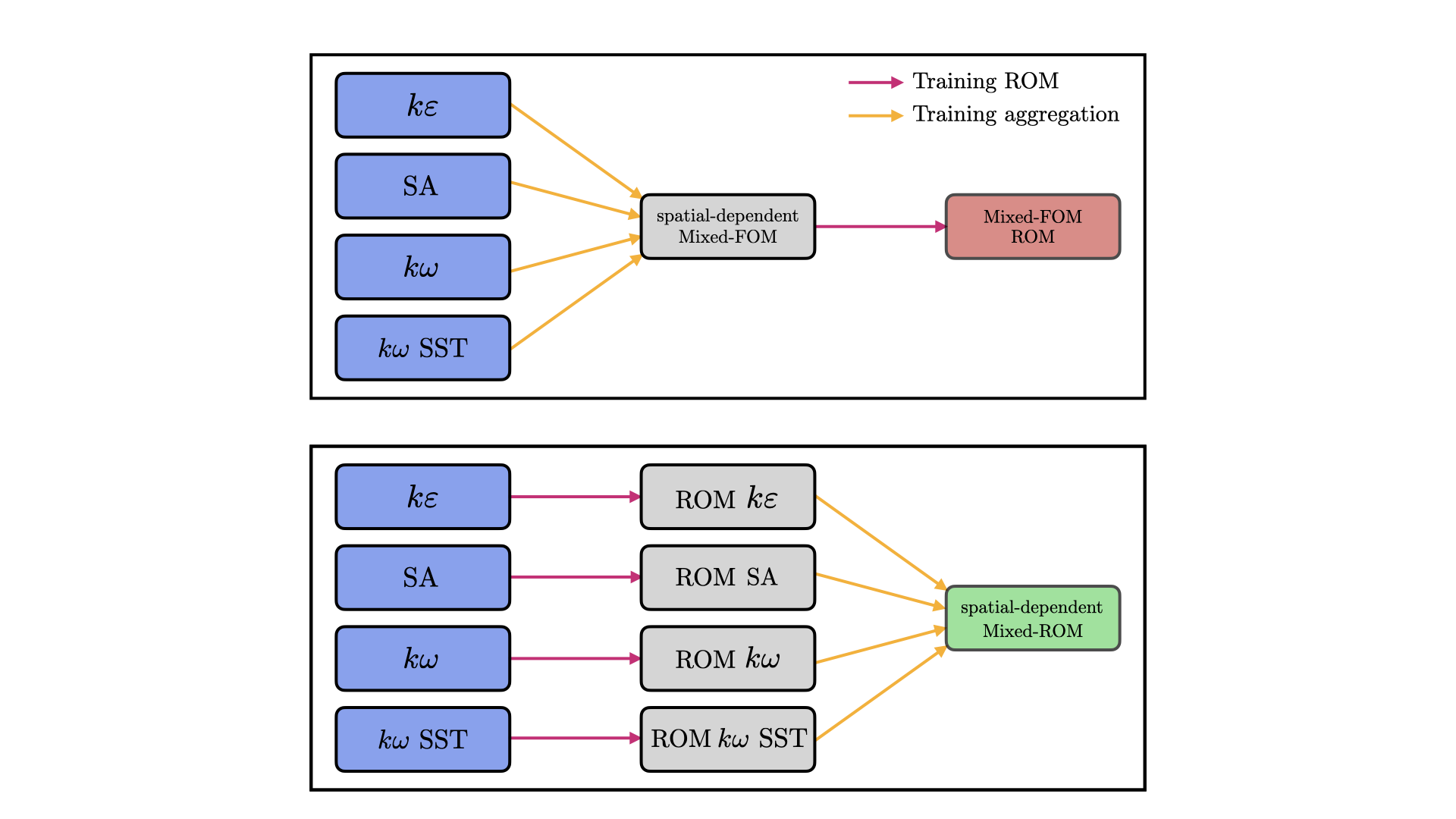}
     \caption{Schematic representation of the mixed aggregation pipelines.}
     \label{pipelines}
\end{figure}

 Since the full-order modeling framework considered in this work involves multiple turbulence models, two main strategies for combining ROMs and spatial aggregation techniques can be adopted. The first strategy consists in applying a space-dependent aggregation directly to the solutions obtained from several RANS simulations, and subsequently training a single ROM on the resulting aggregated high-fidelity fields. Alternatively, one may construct separate ROMs corresponding to each turbulence model and then perform a space-dependent aggregation at the ROM level, combining the reduced-order predictions instead of the full-order solutions. We will refer to the first pipeline as Mixed-FOM ROM (MFR) and to the second as Mixed-ROM (MR). 
A schematic representation of the two overall workflows is shown in Figure \ref{pipelines}.
Both pipelines rely on the selection of $n_M$ eddy-viscosity turbulence models. In this work, four widely used closures in industrial applications are considered: the Spalart-Allmaras  \cite{SA_article}, $k\varepsilon$ \cite{kappaepsilon}, $k\omega$ \cite{kappaomega}, and $k\omega$ SST \cite{kappaomegasst} models. For each turbulence closure, $N_{\mathrm{train}}$ baseline RANS simulations are performed, yielding a total of $N_{\mathrm{train}}^{\mathrm{tot}} = n_M \cdot N_{\mathrm{train}}$ high-fidelity RANS solutions. These full-order simulations constitute the foundation of the proposed framework and are organized into four distinct datasets, each associated with a specific turbulence model and containing $N_{\mathrm{train}}$ snapshots. At this stage of the pipeline, two alternative strategies can be pursued. In the MFR approach, a space-dependent weighting technique is first applied to the four FOM datasets to construct a new set of $N_{\mathrm{train}}$ aggregated high-fidelity solutions. Each aggregated snapshot is obtained as a convex combination of the corresponding RANS solutions, with spatially varying weights optimized to reflect the local predictive accuracy of each turbulence model across the computational domain. A non-intrusive reduced-order model is then trained on the dataset of the aggregated RANS snapshots, yielding a reduced-dimensional representation of a space-dependent aggregated model.

In the MR approach, separate non-intrusive ROMs are instead trained on each of the four RANS datasets, resulting in $n_M = 4$ distinct ROMs, each tailored to reproduce the flow dynamics predicted by a specific turbulence closure. Subsequently, a space-dependent weighting model is trained to learn a mapping from the input features $\boldsymbol{\eta}_{\mathrm{train}}$ to the optimal spatially varying weights $\omega_i(\boldsymbol{\eta}_{\mathrm{train}})$, $i = 1,\dots,n_M$, to be assigned to each ROM component. During the online phase, these ROMs enable rapid predictions for previously unseen configurations, and their outputs are combined through the learned weighting strategy. By exploiting the complementary strengths of the different turbulence models, this approach is expected to improve predictive accuracy compared to any single ROM.

In the online phase of the MFR pipeline, the input parameters are directly provided to the final aggregated ROM, without the need to execute multiple RANS simulations. Similarly, in the MR pipeline, the pre-trained ROMs are evaluated for unseen configurations, and their outputs are aggregated using the previously trained KNN or ANN regressors. In both cases, the final prediction relies exclusively on pre-trained models and does not require additional high-fidelity simulations, resulting in substantial computational speed-ups.

Before presenting the computational results, it is useful to discuss the expected differences between the two approaches. On the one hand, the MFR approach relies on training a single ROM on space-dependently aggregated high-fidelity data. This choice naturally reduces the offline computational cost associated with ROM training. Although this cost is not particularly limiting in the present work, it may become significant as the dimensionality of the solution space or the complexity of the physical model increases. On the other hand, the MR approach requires the training of multiple ROMs, one for each turbulence model, which is inherently more computationally demanding in the offline phase. However, since the resulting model mixture is explicitly designed to optimally combine the strengths of all component ROMs, the final aggregated predictions are expected to exhibit improved accuracy compared to those obtained from any single RANS model or individual ROM. Indeed, by leveraging the complementary characteristics of different turbulence closures, the mixed-ROM framework provides a more robust and accurate representation of the underlying flow physics. 

Finally, note that both pipelines rely on a combination of space-dependent aggregation learning, trained on high-fidelity DNS/LES data, and model order reduction, built from RANS simulations. The former typically requires only a limited amount of high-fidelity data, whereas the latter demands a substantially larger dataset, although at lower fidelity. From this perspective, the proposed framework can be interpreted as an efficient strategy for exploiting multiple levels of fidelity in the simulation of turbulent flows.

\section{Numerical results}
\label{sec:results}
\par In this Section, we present the results of the numerical simulations conducted to assess the performance of the surrogate model mixtures derived from the two aggregation pipelines we propose. The test cases considered are the two-dimensional, incompressible flow over periodic hills, and the flow over a bump, which are well-established benchmark problems in CFD \cite{periodic_hill_extra1,periodic_hill_extra2,periodic_hill_extra3,Xiao,Breuer,Frohlich,Rapp,CinnellaLES,matai2019large,mcconkey2021curated,wu2023enhancing}. We select as physical quantity of interest the horizontal component $U_x$ of the two-dimensional velocity field. The FOMs are developed using OpenFOAM \cite{Jasak}, an open-source \texttt{C++} CFD software based on finite-volume discretization \cite{Finite_volume}.
All ROMs are constructed within the \texttt{EZyRB} \cite{ezyrb} framework using an AE for dimensionality reduction, coupled with RBF interpolation as approximation method.

\subsection{Test case 1: periodic hills}
\label{sec:hills-case}

\par We consider a two-dimensional, incompressible flow over periodic hills with parameterized geometries. This test case is well documented in the CFD literature, as it is supported by extensive experimental and numerical high-fidelity data across a wide range of geometrical parameters and regimes, from $\mathrm{Re}=700$ to $\mathrm{Re}=37000$ \cite{Xiao,Breuer,Frohlich,Rapp,CinnellaLES}. In this study, we fix the Reynolds number to $\mathrm{Re}=5600$, which is a classical choice due to the abundance of available reference data \cite{Xiao,Breuer,Rapp}.
\par
\begin{figure}[!htpb]
     \centering
     \includegraphics[width=0.7\textwidth]{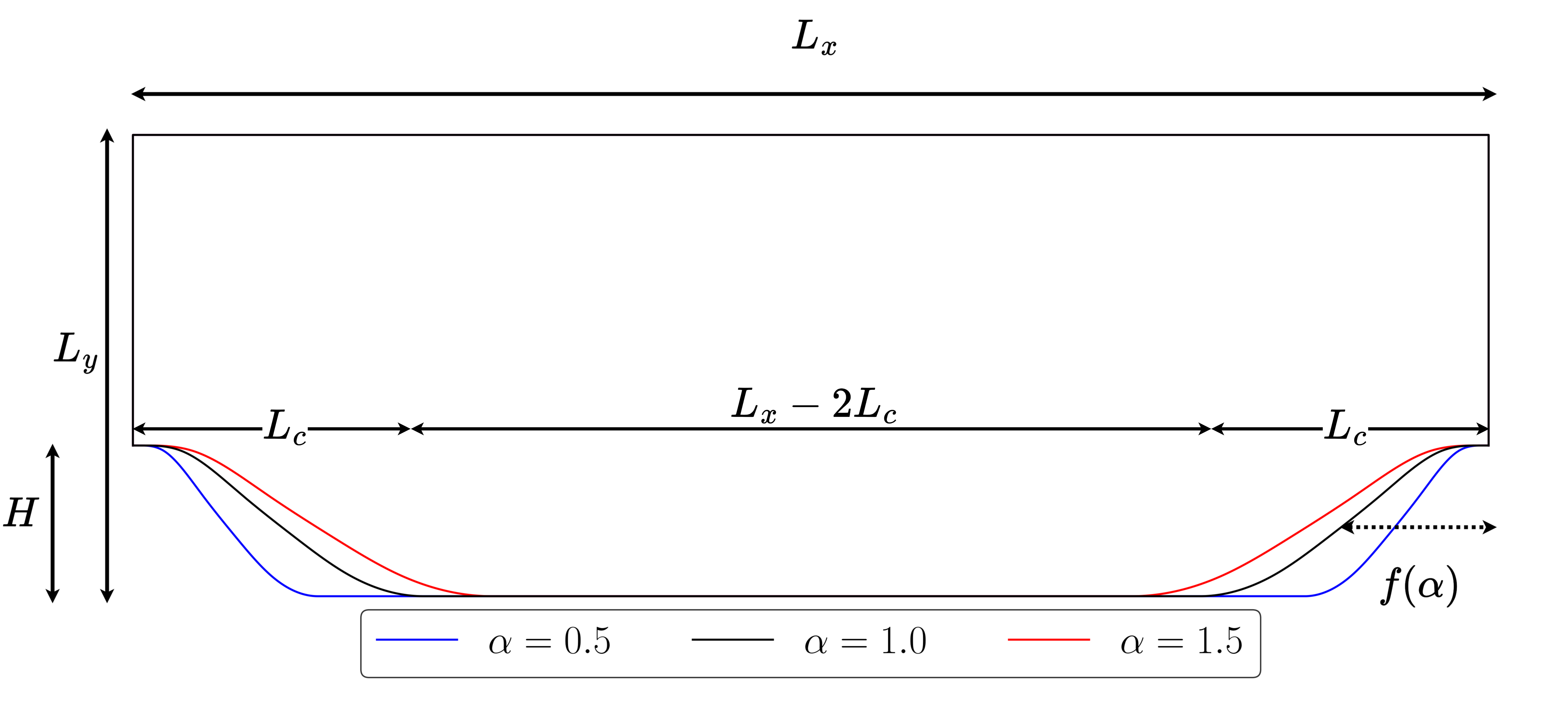}
     \caption{\emph{Test case 1}. Schematic representation of the geometrical parametrization of the computational domain. Hill profiles corresponding to three different values of the hill stretch factor $\alpha$ are illustrated. The quantity $L_c$ indicates the streamwise extent of the curved regions of the hill profile for the case $\alpha=1$.}
     \label{periodic_hill}
\end{figure}

 The geometrical parametrization of the computational domain involves four features: the crest height $H$, the domain height $L_y$, the domain length $L_x$ and the hill stretch factor $\alpha$. All geometric dimensions are normalized with respect to the crest height $H$, with $L_y/H=3.036$ and $H=1\ \mathrm{m}$ fixed across all configurations. The $x-$ and $y-$coordinates are aligned with the streamwise and wall-normal directions, respectively. Variations in $\alpha$ modify the hill shape through a piecewise cubic spline, producing progressively smoother and more elongated profiles as $\alpha$ increases (see Figure \ref{periodic_hill}).
\par A total of $N_{\mathrm{tot}}=11$ geometrical configurations is generated by varying the physical parameters $\boldsymbol{\mu}=[\alpha,L_x]$. The \textit{training} set comprises $N_{\mathrm{train}}=9$ configurations corresponding to $\alpha=0.5,1,1.5$, each combined with three streamwise lengths defined by the following linear relation:
\begin{equation*}
    \displaystyle{\frac{L_x}{H}=3.858\alpha+c},\quad\text{with}\quad c=[2.142,5.142,8.142].
\end{equation*}
\par $N_{\mathrm{test}}=2$ additional intermediate cases with $\alpha=0.75,1.25$ and the central value of $\displaystyle{\frac{L_x}{H}}=3.858\alpha+5.142$ constitute the \textit{test} set.
\par All these geometries exactly match the configurations of the $N_{\mathrm{tot}}$ DNS solutions of Xiao et al. \cite{Xiao} at $\mathrm{Re}=5600$, which we adopt as reference data.
\par The problem instances in the training set are employed in the offline phases of both aggregation pipelines to run the $N_{\mathrm{train}}^{\mathrm{tot}}$ RANS simulations, train all reduced order models and tune the space-dependent weighting techniques using $N_{\mathrm{train}}$ DNS snapshots. In contrast, the test configurations are used in the online phases to evaluate the predictive performance of individual models and aggregated mixtures on unseen parameter values and to compare them against the corresponding $N_{\mathrm{test}}$ DNS solutions.

\subsubsection{RANS FOMs}
\par All computational meshes are obtained from a baseline configuration of the domain with $\alpha=1$ and $L_x/H=9$, by applying a linear deformation driven by the selected parameters $\boldsymbol{\mu}$. The hill profile is partitioned into two curved segments and one flat region (see Figure \ref{periodic_hill}); the curved portions are scaled proportionally to $\alpha$. Streamwise coordinates are redistributed linearly within each segment, while the wall-normal coordinates are computed from the cubic spline definition of the hill.
\par Periodic boundary conditions are imposed at the inlet and outlet boundaries. No-slip Dirichlet boundary conditions are enforced for the velocity field on the top and bottom walls, while homogeneous Neumann boundary conditions are applied for the pressure on the same boundaries. The coupling between the pressure and velocity fields is handled using the Semi-Implicit Method for Pressure-Linked Equations (SIMPLE) algorithm \cite{SIMPLE} for steady flows.
\par
\begin{figure}[!htpb]
     \centering
     \includegraphics[width=0.6\linewidth]{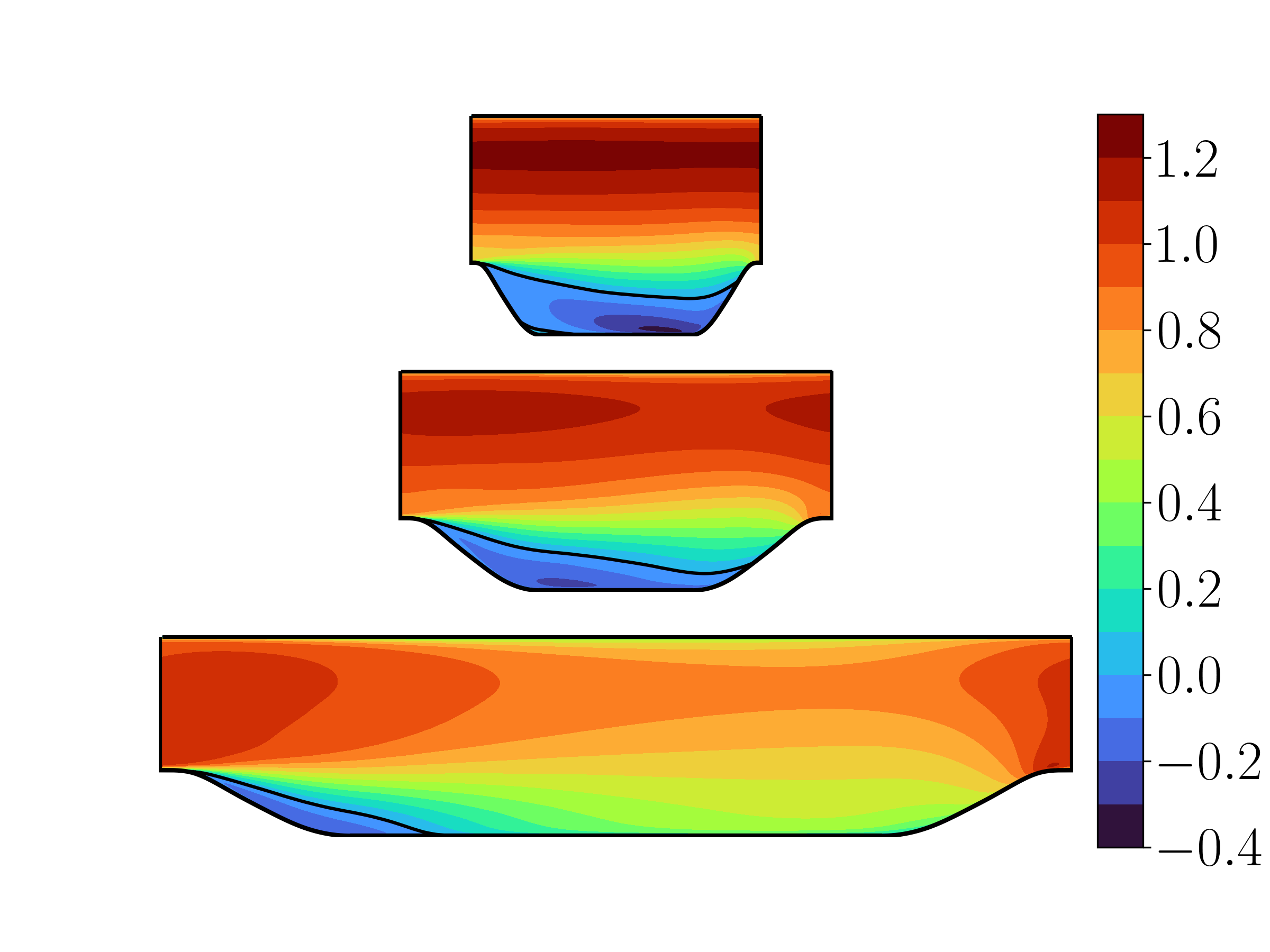}
     \caption{\emph{Test case 1}. Deformed representations of the periodic hill configurations with $\alpha=0.5$ and $L_x/H=7.071$ (top), $\alpha=1.0$ and $L_x/H=9.0$ (middle) and with $\alpha=1.0$ and $L_x/H=12.0$ (bottom).}
     \label{two_meshes}
\end{figure}

 Each mesh is composed of $N_h=32000$ cells, with grid refinement applied towards both the top and bottom walls. Importantly, all RANS simulations share the same number of degrees of freedom to ensure consistency in the subsequent construction of reduced order models based on the corresponding solution snapshots. However, the spatial coordinates $\boldsymbol{x}$ vary across the simulations, reflecting the distinct geometrical configurations associated with each mesh (see Figure \ref{two_meshes}). Finally, the DNS reference solutions are linearly interpolated onto the coarser grids used for the RANS simulations. This pre-processing step enables a consistent pointwise comparison across datasets and allows the training of the space-dependent weighting strategies on a unified spatial grid.

\par Each problem instance in the training set is simulated using $n_M=4$ different RANS turbulence models, which constitute the baseline for the subsequent aggregation procedures. Specifically, the \textit{Spalart-Allmaras} \cite{SA_article}, $k\varepsilon$ \cite{kappaepsilon}, $k\omega$ \cite{kappaomega} and $k\omega\ \text{SST}$ \cite{kappaomegasst} models are selected, as they are among the most widely adopted RANS closures in industrial CFD applications. With $N_{\mathrm{train}}=9$ simulations per model, this results in $N_{\mathrm{train}}^{\mathrm{tot}}=36$ high-fidelity RANS simulations, from which the corresponding full-order solutions are extracted. These snapshots form the four distinct high-fidelity datasets, one for each turbulence model, that constitute the foundation for both aggregation pipelines.
\begin{figure}[!htpb]
     \centering
     \includegraphics[width=0.8\textwidth]{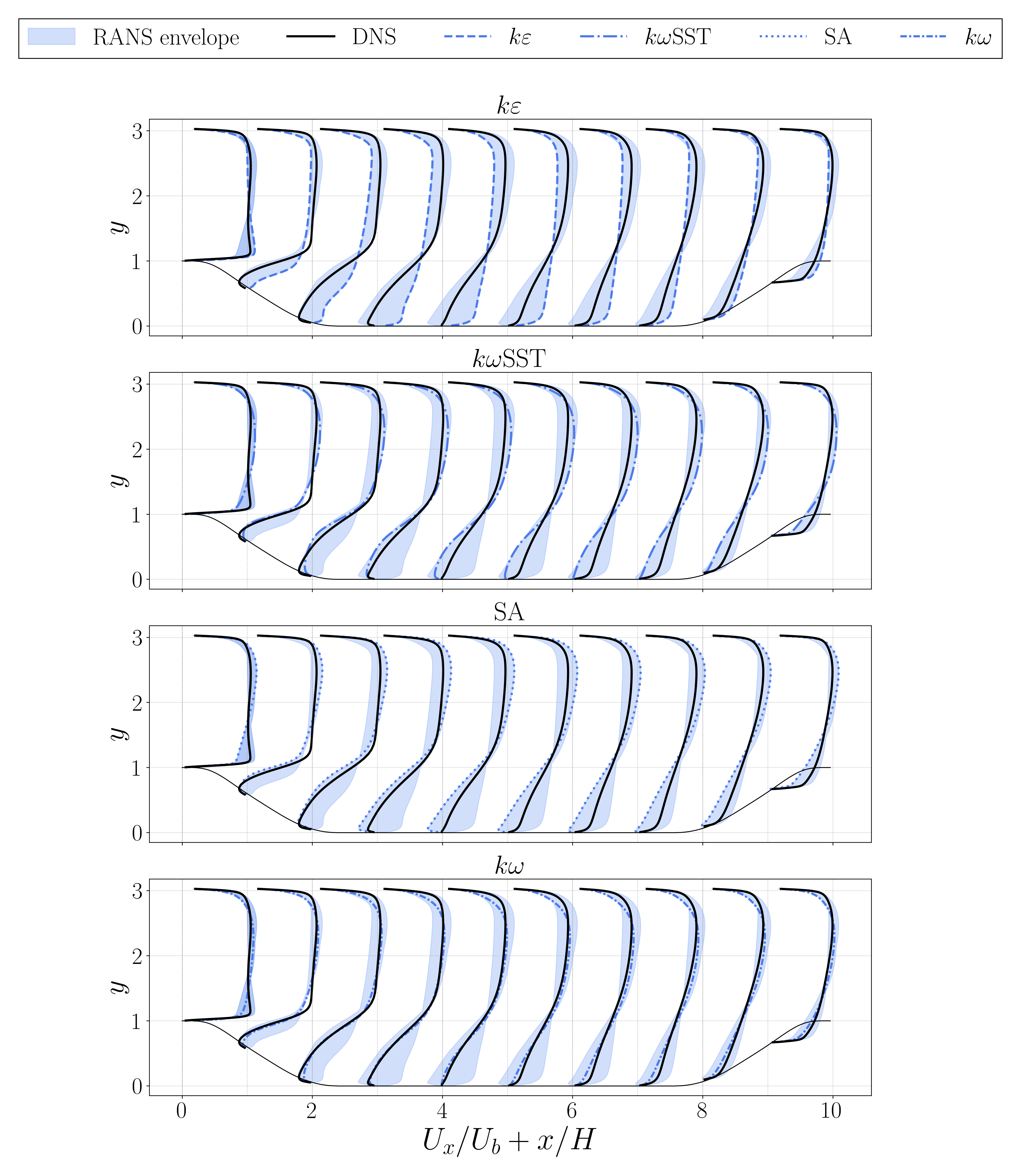}
     \caption{\emph{Test case 1}. Comparison of velocity profiles obtained from RANS simulations using different turbulence models and DNS reference data. We display the horizontal component of the velocity field, extracted at ten distinct downstream locations. The results correspond to the test configuration with $\alpha=1.25$ and $L_x/H=9.96$. The light blue area indicates the envelope of all the different RANS simulations.}
     \label{Ux_RANS_comp}
\end{figure}
\par Figure \ref{Ux_RANS_comp} compares the horizontal velocity profiles from the four RANS models against the corresponding DNS data for one of the test configurations. The four RANS turbulence models yield significantly different predictions. Moreover, the accuracy of each closure varies considerably across specific regions of the computational domain. Notably, this spatial variability in performance naturally motivates the concept of locally selecting the best-performing model throughout the domain. Furthermore, the DNS solution lies almost entirely within the envelope spanned by the four RANS predictions, suggesting that convex combinations of the individual model outputs can exploit their complementary strengths and yield improved accuracy.

\subsubsection{Aggregated models}\label{sec:aggregation-res}
\par This Subsection presents the results obtained from the aggregated models produced by the two pipelines we proposed. We evaluate model performance on the test set in terms of predictive accuracy, thereby measuring the ability of the surrogate mixtures to effectively generalize to previously unseen flow configurations.
Regarding the KNN regression, we set the number of neighbors for training to $K=4$, after hyperparameter tuning.
The predictions are compared against both DNS data and the baseline RANS models, enabling a clear assessment of the improvements achieved through the aggregation framework.
\par
\begin{figure}[!htpb]
     \centering
     \includegraphics[width=0.75\textwidth]{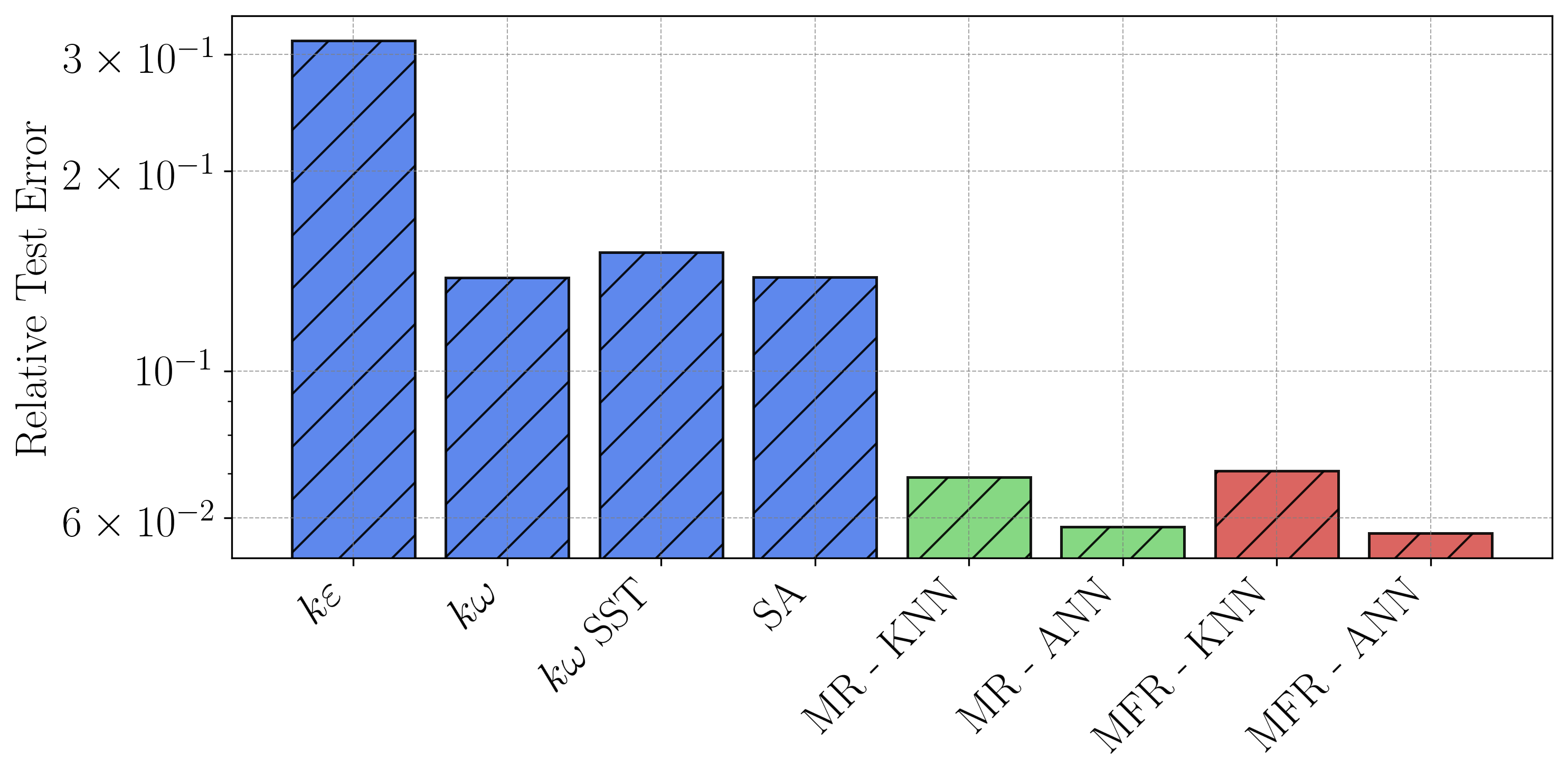}
     \caption{\emph{Test case 1}. Comparison of relative test errors with respect to DNS data for baseline RANS models and surrogate model mixtures.}
        \label{Relative_errors}
\end{figure}

 As a first quantitative comparison, Figure \ref{Relative_errors} reports the test errors of the individual baseline RANS models and of all the aggregated surrogate models we developed, computed relative to the DNS reference solutions. The metric adopted is the mean relative error in the $L^2$ norm. It is evident that all aggregated models, independently of the chosen weighting strategy or the specific aggregation pipeline, consistently and substantially outperform the baseline RANS component models.
\par Additional observations can be drawn from the comparison in Figure \ref{Relative_errors}. First, as we expected, the surrogate models obtained with the ANN weighting strategy outperform those based on the KNN supervised regression technique. The superiority of the ANN-based weighting is consistently observed across both aggregation pipelines. Second, for a fixed weighting strategy, the surrogate models produced by the two different pipelines exhibit broadly comparable predictive accuracies. This result indicates that, in practice, the two pipelines are essentially equivalent in terms of the quality of the surrogate mixtures they generate. However, while in MFR we only have to train a single AE, in MR we have to train $n_M=4$ AEs, thus multiplying by four the ROM training time. In the remaining of the paper, therefore, we will only consider the MFR approach, which, involves less demanding trainings, while keeping a comparable accuracy to MR.


\par As an additional validation metric commonly adopted in the literature \cite{periodic_hill_extra3}, we consider the reattachment point $x_R$. Table \ref{reattachment_points} reports the values of $x_R$ extracted from the velocity fields predicted by the individual RANS turbulence models, together with those obtained from the two surrogate model mixtures constructed using the ANN weighting strategy. For completeness, the Table also includes the relative $L^2$ error of each prediction with respect to the DNS reference. In particular, the results correspond to the test case characterized by $\alpha=1.25$ and $L_x/H=9.96$, hereafter denoted as \texttt{T1}. For this configuration, the DNS value of the reattachment point is $x_R^{DNS}/H=4.33$.
\begin{table}[!htpb]
    \caption{\emph{Test case 1}. Comparison of reattachment point predictions from DNS, baseline RANS models and surrogate models constructed via ANN weighting. The corresponding relative errors with respect to the DNS reference are also reported. The results refer to the \texttt{T1} test configuration.}
    \centering
    \begin{tabular}{c c c}
    \toprule
         \makecell[c]{\textbf{Model}} & \makecell[c]{$\boldsymbol{x_R/H}$} & \makecell[c]{\textbf{Relative error}} \\
         \midrule
    $\text{DNS}$ & $4.33$ & --- \\
    \midrule
    $k\varepsilon$ & --- & --- \\
    $k\omega$ & $4.25$ & $0.019$\\
    $k\omega\ \text{SST}$ & $5.82$ & $0.344$\\
    $\text{SA}$ & $7.23$ & $0.670$\\
    MFR - ANN & $\boldsymbol{4.40}$ & $\boldsymbol{0.016}$\\
    MR - ANN & $4.41$ & $0.019$\\
    \bottomrule
    \end{tabular}
    \label{reattachment_points}
\end{table}
\par It is clear that the $k\varepsilon$ model completely fails to capture the presence of a recirculation region, as reflected by the absence of any identifiable reattachment point. Among the baseline turbulence closures, only the $k\omega$ model provides a prediction that closely matches the DNS reference value. This behavior is fully consistent with expectations, given the ability of the $k\omega$ formulation to deliver superior performance close to the walls and in separated flow regions \cite{RANS_REF1,ref_turbulence_models,kappaomegaref}. More importantly, the aggregated surrogate models not only reach but, in the case of the former, slightly surpass the accuracy of the best-performing individual RANS closure. Moreover, despite being generated through different aggregation pipelines, the two surrogate mixtures yield predictions that are again remarkably similar, confirming the robustness of the aggregation framework. Overall, this complementary analysis reinforces the conclusion that the final surrogate mixtures successfully exploit the complementary strengths of the underlying turbulence models across the entire domain. As a result, they deliver enhanced predictive accuracy not only for previously unseen configurations, but also for flow quantities that were not directly employed during their training phase.

\begin{figure}[!htpb]
     \centering
     \includegraphics[width=0.75\textwidth]{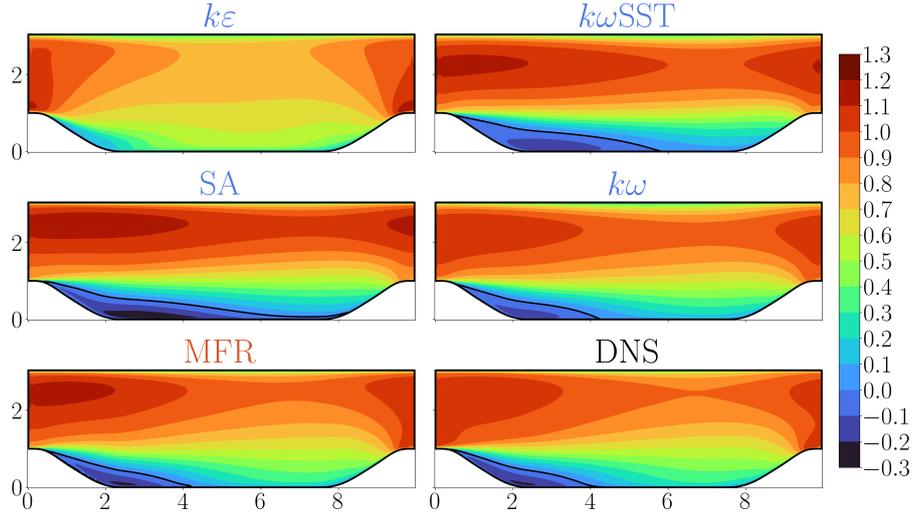}
     \caption{\emph{Test case 1}. Comparison of the horizontal velocity component predicted by the individual RANS models and the MFR-ANN surrogate, with DNS reference solution. The results refer to the \texttt{T1} test configuration.}
        \label{preds_MFROM}
\end{figure}

\par
\begin{figure}[!htpb]
     \centering
     \includegraphics[width=0.75\textwidth]{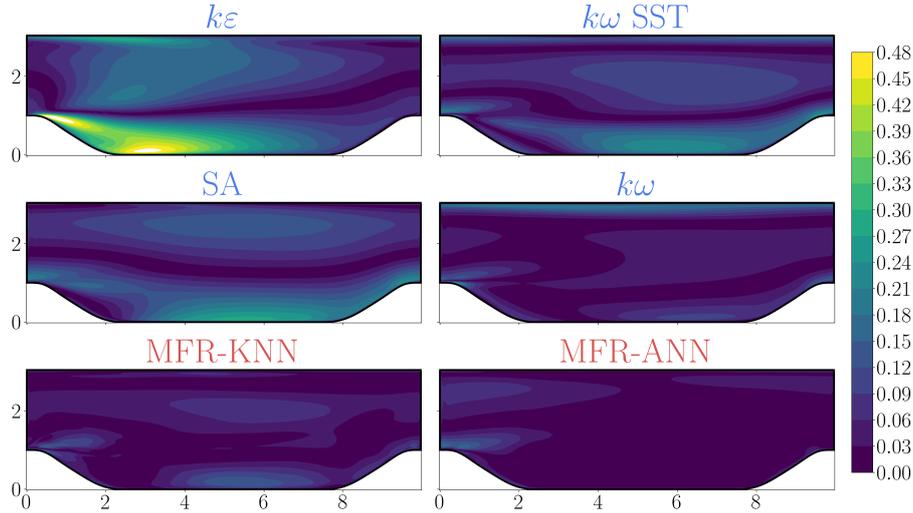}
     \caption{\emph{Test case 1}. Absolute error maps with respect to the DNS benchmark for the horizontal velocity component, predicted by the individual RANS models, the MFR-KNN and the MFR-ANN surrogates. The results correspond to the \texttt{T1} test configuration.}
        \label{ABS_ERRORS_MFROM}
\end{figure}

 Figure \ref{preds_MFROM} compares the horizontal velocity fields predicted by the individual RANS models, the MFR-ANN surrogate and the DNS. In particular, the results refer to the \texttt{T1} test configuration. In addition, Figure \ref{ABS_ERRORS_MFROM} presents the corresponding absolute error maps. These errors are computed with respect to the DNS solution for the individual RANS models as well as for the MFR–KNN and MFR–ANN surrogates.
\par The MFR prediction appears overall the closest to the reference solution. Moreover, the surrogate models not only achieve visibly improved accuracy, but also substantially reduce the errors across the entire computational domain. Notably, it is evident that the KNN-based mixture, while still significantly improving the predictive performance upon all individual RANS component models, generally exhibits larger errors than its ANN-based counterpart.

\par We highlight that the MFR-ANN is never exposed to the individual RANS solutions shown in the plots above. In fact, these simulations are independently run solely to enable a direct comparison between the aggregated model and the standalone turbulence models. The surrogate is instead trained on high-fidelity aggregated predictions, obtained as convex combinations of RANS outputs computed for different geometrical configurations of the domain. This ensures that the improved predictive accuracy we observe is not the result of the ROM learning the behavior of any specific turbulence model in the \texttt{T1} configuration, but rather of its ability to generalize from aggregated information across multiple, distinct training geometries.
\par
\begin{figure}[!htpb]
   \centering
   \begin{subfigure}[b]{0.3\linewidth}
       \centering
       \includegraphics[width=\linewidth]{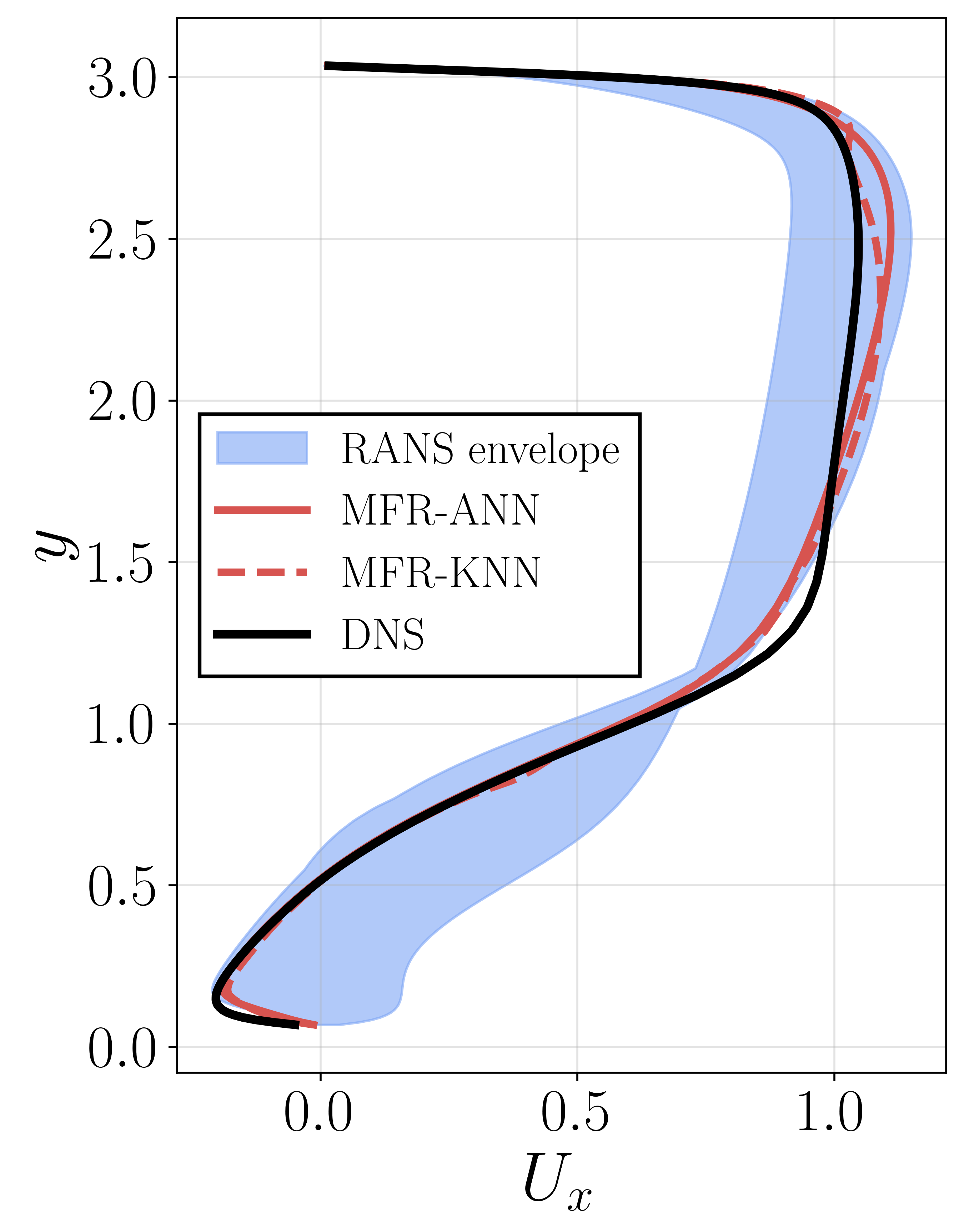}
       \caption{$U_x$ velocity profiles at $x/H=2$.}
       \label{profiles_MFROM_U_2}
   \end{subfigure}
   \hfill
    \begin{subfigure}[b]{0.3\linewidth}
        \centering        
        \includegraphics[width=\linewidth]{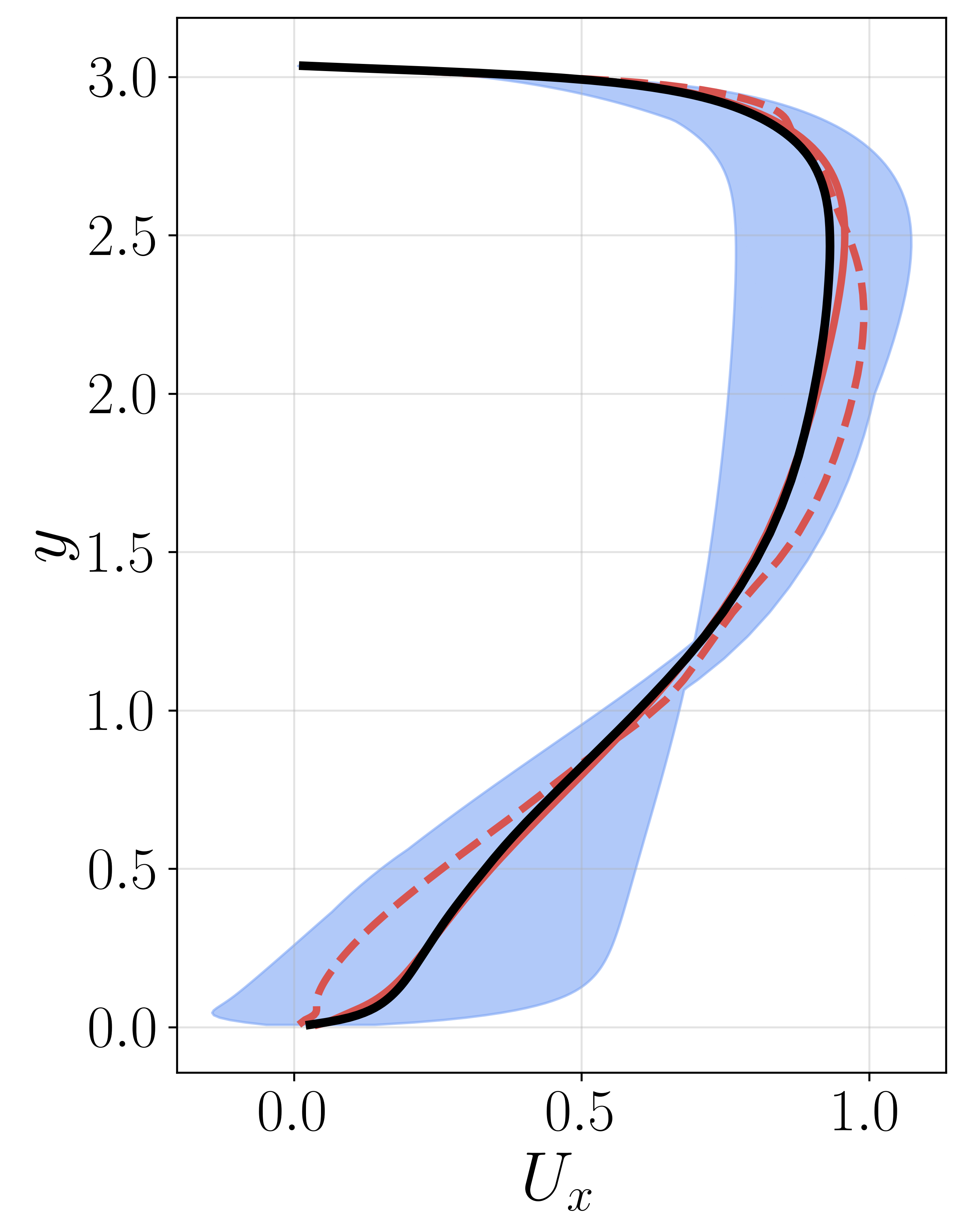}
        \caption{$U_x$ velocity profiles at $x/H=5$.}
        \label{profiles_MFROM_U_5}
    \end{subfigure}
   \hfill
    \begin{subfigure}[b]{0.3\linewidth}
        \centering        
        \includegraphics[width=\linewidth]{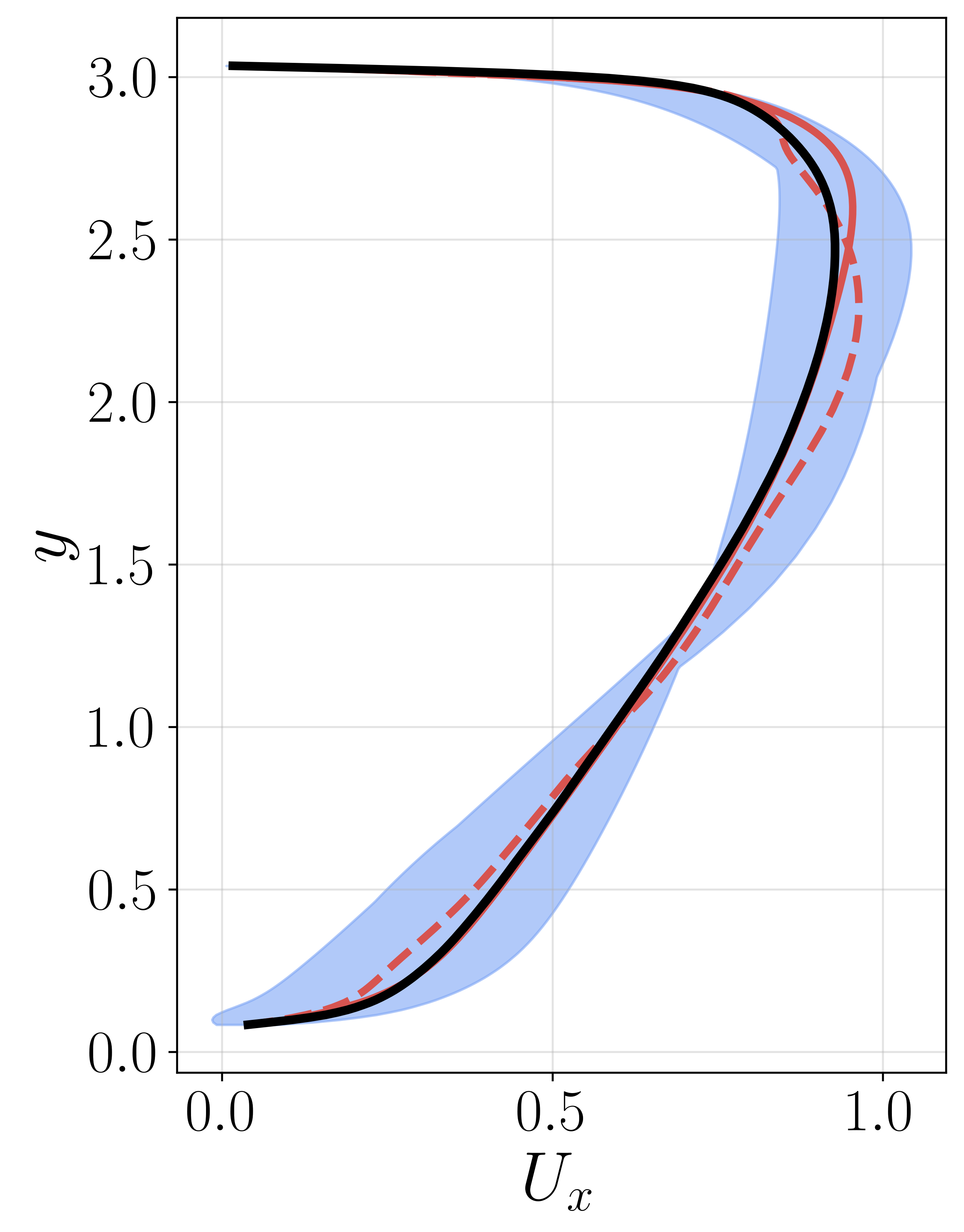}
        \caption{$U_x$ velocity profiles at $x/H=8$.}
        \label{profiles_MFROM_U_8}
    \end{subfigure}
    \caption{\emph{Test case 1}. Comparison of horizontal velocity profiles obtained from the two MFR models employing different weighting techniques, alongside DNS reference data. The profiles are extracted at $x/H=2$ (left), $x/H=5$ (center) and $x/H=8$ (right) for the \texttt{T1} test configuration. The shaded region indicates the prediction range of the four individual RANS turbulence models.}
    \label{profiles_MFROM_U}
\end{figure}

 Figure \ref{profiles_MFROM_U} presents a comparison of the horizontal velocity profiles obtained from the MFR models developed using the KNN and ANN weighting strategies, alongside the associated DNS solution. Specifically, the profiles are extracted at three distinct downstream locations: $x/H=2$, located within the recirculation region, $x/H=5$, near the center of the domain, and $x/H=8$. The shaded area in each plot represents the range spanned by the predictions of the four individual RANS turbulence models for this specific problem instance, providing a visual indication of the improvement achieved by the aggregated surrogates.

\par Figure \ref{profiles_MFROM_U} illustrates that both surrogate models produce smooth solutions that are relatively close to one another. Importantly, the predicted velocity profiles exhibit relatively strong agreement with the DNS data, effectively capturing the reference behavior despite the overall considerable variability observed among the individual RANS component models, as indicated by the wide shaded region. Notably, the ANN-based approach provides the closest match to the reference, whereas the KNN-based model tends to show the largest deviations, visually reinforcing the findings of the previous analysis.
\par
\begin{figure}[!htpb]
    \centering
    \begin{subfigure}[b]{\linewidth}
        \centering        
        \includegraphics[width=0.75\textwidth]{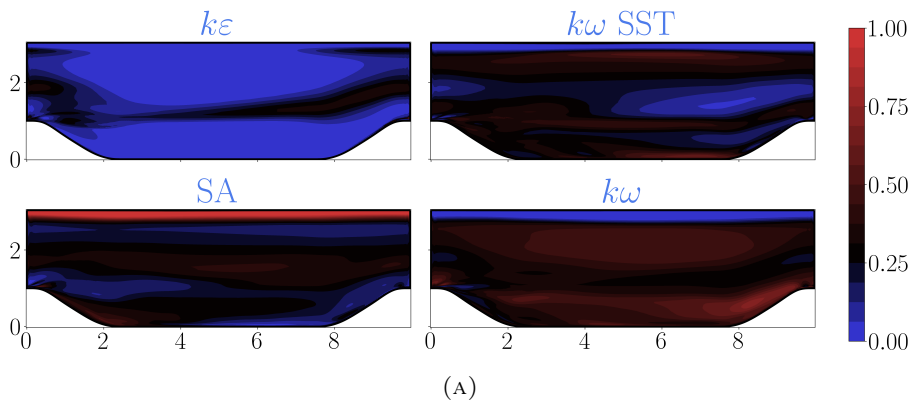}
        \caption{}
        \label{weightsKNN_U}
    \end{subfigure}
    \caption{\emph{Test case 1}. Spatial distribution of the weighting functions for the four component RANS models, in the MFR-KNN pipeline.}
    \label{weightsKNN_MF}
\end{figure}

 Figure \ref{weightsKNN_MF} shows the spatial distributions of the weighting functions assigned to the four component RANS turbulence models, as determined by the MFR-KNN approach.

\par The space-dependent weights derived from the KNN weighting strategy show that the SA model receives very high weights in the boundary layer near the upper wall; the $k\omega$ model is strongly activated in the lower portion of the domain, while the $k\varepsilon$ model is consistently assigned very low weights, reflecting its generally poorer performance. Overall, the weight distributions align with the local predictive capabilities of each turbulence closure, as the KNN weighting technique effectively identifies the best-performing RANS model in each subregion of the computational domain.
\par
\begin{figure}[!htpb]
     \centering
     \includegraphics[width=0.75\textwidth]{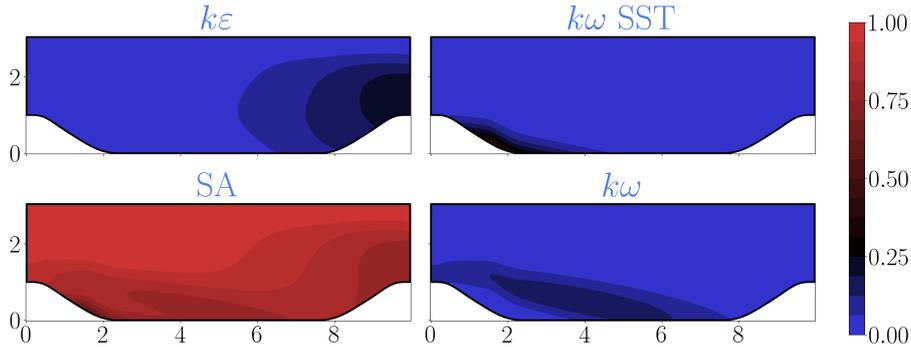}
     \caption{\emph{Test case 1}. Spatial distributions of the weighting functions for the four component RANS models, in the MFR-ANN pipeline.}
        \label{weightsANN_MF}
\end{figure}

 Finally, Figure \ref{weightsANN_MF} depicts the spatially varying weights assigned to the four individual RANS turbulence models, as predicted by the ANN weighting approach. It can be noticed that in this case the optimal spatial distributions of the weighting functions exhibit substantially less variability compared to those produced by the KNN approach. The ANN-based method tends to assign very high weights to a single best-performing model over large portions of the domain, while the remaining closures contribute only marginally or are entirely disregarded. As a result, it places less emphasis on capturing model uncertainty and the weighting patterns are less interpretable.

\begin{remark} We remark that, within this pipeline, the space-dependent weights are used only to build the high-fidelity aggregated predictions that constitute the training dataset (Mixed-FOM) of the non-intrusive reduced-order model. During the online phase, the surrogates do not actually employ these weight fields, which are provided exclusively to show the effectiveness of the weighting strategies.
\end{remark}

\subsection{Test case 2: flow over a bump}
\label{sec:bump-case}
As a second test case to validate our methodology, we consider the LES dataset for flow over a parametric family of surface bumps by Matai and Durbin~\cite{matai2019large}, recently made available for data-driven turbulence closure development. The geometry consists of a circular-arc bump with convex fillets, parameterized by its height $h$. The Reynolds number based on inlet momentum thickness is fixed at $Re_{\theta}=\num{2500}$, while the bump-height Reynolds number $Re_h$ varies from approximately $\num{13250}$ to $\num{27850}$.

For $h=20\,\mathrm{mm}$ the flow remains attached, whereas at $h=26\,\mathrm{mm}$ mild separation occurs, and at $h=42\,\mathrm{mm}$ a small separated region forms downstream of the bump. Unlike massively separated benchmark cases (e.g., periodic hills), these configurations exhibit mild separation characterized by a high turbulent kinetic energy region that lifts away from the wall upstream of separation, attributed to the adverse pressure gradient and the resulting mean shear. The dataset also captures strong non-equilibrium turbulence effects. Further details regarding the simulations set up can be found in \cite{mcconkey2021curated}.

The DNS dataset includes five different bump heights: $h=20, 26, 31, 38, 42\, \mathrm{mm}$.

We call $H$ the global domain height, which is fixed to $H=150\, \mathrm{mm}$ in all the simulations. In all the results, we show the normalized geometry with respect to $H$.
We consider a geometrical parametrization with one parameter $h$, namely the bump height.

\par In addition to the different separation characteristics between the periodic hills and the parametrized bumps test cases, there are also significant differences in terms of data availability. In the present case, the amount of high-fidelity training data is considerably smaller, as we only have $5$ DNS available solutions. As previously discussed, the limited availability of high-fidelity data primarily affects the spatial aggregation procedure, whereas the ROM can, in principle, be trained using a much larger number of snapshots.
To keep an acceptable ROM accuracy, we build the ROM (or ROMs, depending on the chosen pipeline) starting from $N_{\mathrm{train}}^{(mix)}=10$ snapshots for each turbulence model.
The set of $5$ DNS is instead divided into a set of $N_{\mathrm{train}}^{(mix)}=3$ solutions, used to train the aggregation, and $N_{\mathrm{test}}=2$ solutions, used to test the performance of our approach.
In particular, we consider $h=20, 31, 42 \, \mathrm{mm}$ as parameters to train the aggregation, and $h=26, 38 \, \mathrm{mm}$ for testing.

The following part is dedicated to the presentation of the FOM setup for the $N_{\mathrm{train}}^{(mix)}$ snapshots.

\subsubsection{RANS FOMs}
\begin{figure}[!htpb]
     \centering
     \includegraphics[width=0.7\textwidth]{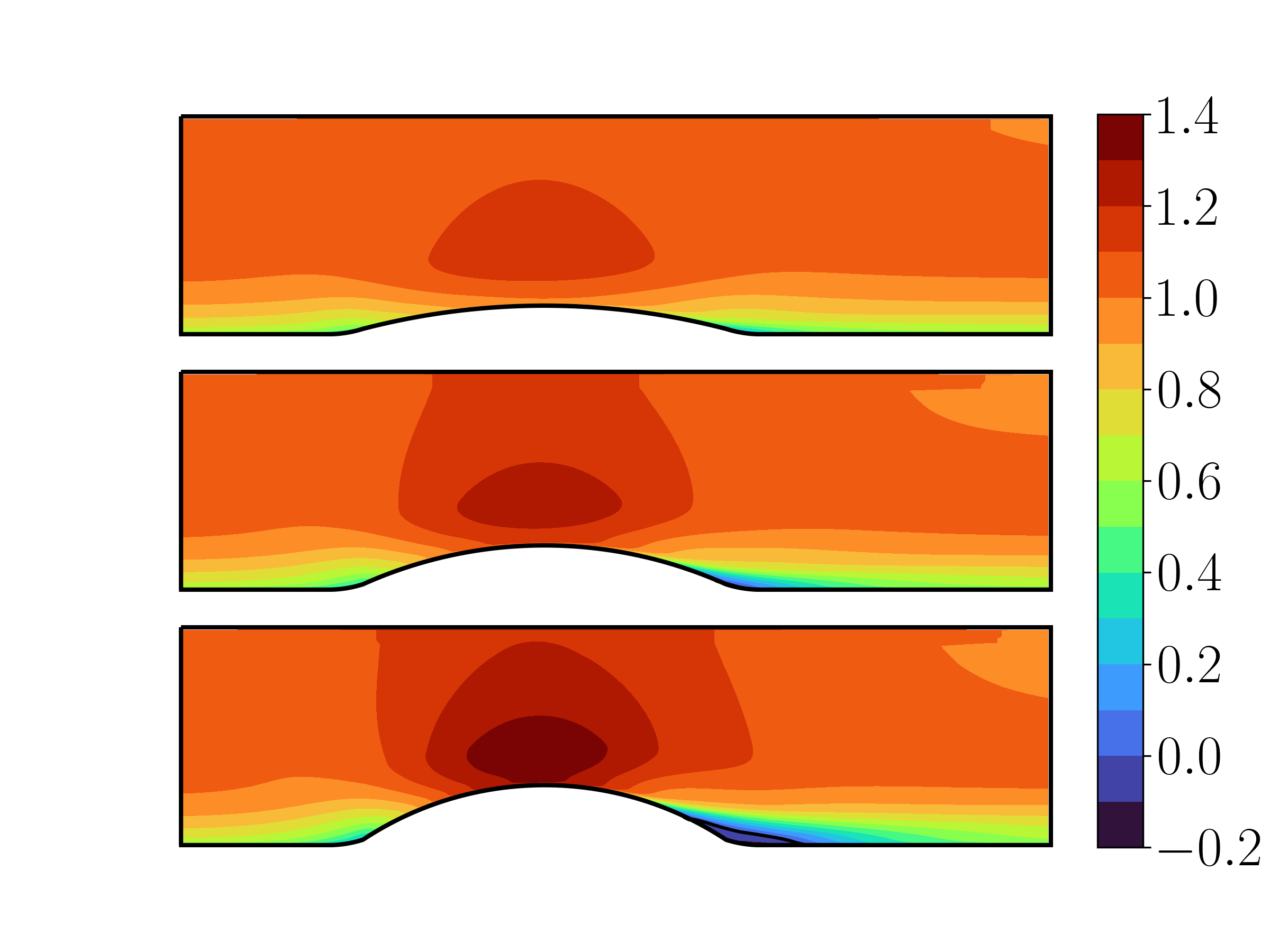}
     \caption{\emph{Test case 2}. Deformed representations of the bumps for different bump heights, namely $h=20, 31, 42 \, \mathrm{mm}$.}
     \label{ref-DNS-bumps}
\end{figure}

Figure \ref{ref-DNS-bumps} shows the reference DNS $U_x$ solution for three different height values. The Figure highlights the differences in the flow behavior, with no separation at $h=20\, \mathrm{mm}$ and a small separated region close to the bump at $h=42 \, \mathrm{mm}$.

\par In this second test case, we employ the same set of RANS turbulence models adopted in the previous case, namely the \textit{Spalart--Allmaras}, $k\varepsilon$, $k\omega$, and $k\omega$~SST closures.

\par As in the first test case, we perform a set of steady-state simulations, where the horizontal velocity field is initialized with a fixed value $U_b=16.683\, \mathrm{m}/\mathrm{s}$, while we consider a zero-valued initial pressure field. We consider a Neumann pressure condition for all the boundaries, except for the outlet, where we impose $p=0$.
Regarding the velocity, we consider a no-slip boundary condition on the bottom, a Neumann condition on the top and on the bottom boundaries. Since the original DNS simulations have a free-stream velocity $U_{\infty}=16.77 \, \mathrm{m}/\mathrm{s}$, we prescribe the inlet velocity field to match this reference value, as done in~\cite{mcconkey2021curated}.
As previously mentioned, we consider $10$ simulations, for $h$ equispaced values within the range $[13, 49] \, \mathrm{mm}$.

\begin{figure}[!htpb]
     \centering
     \includegraphics[width=0.75\textwidth]{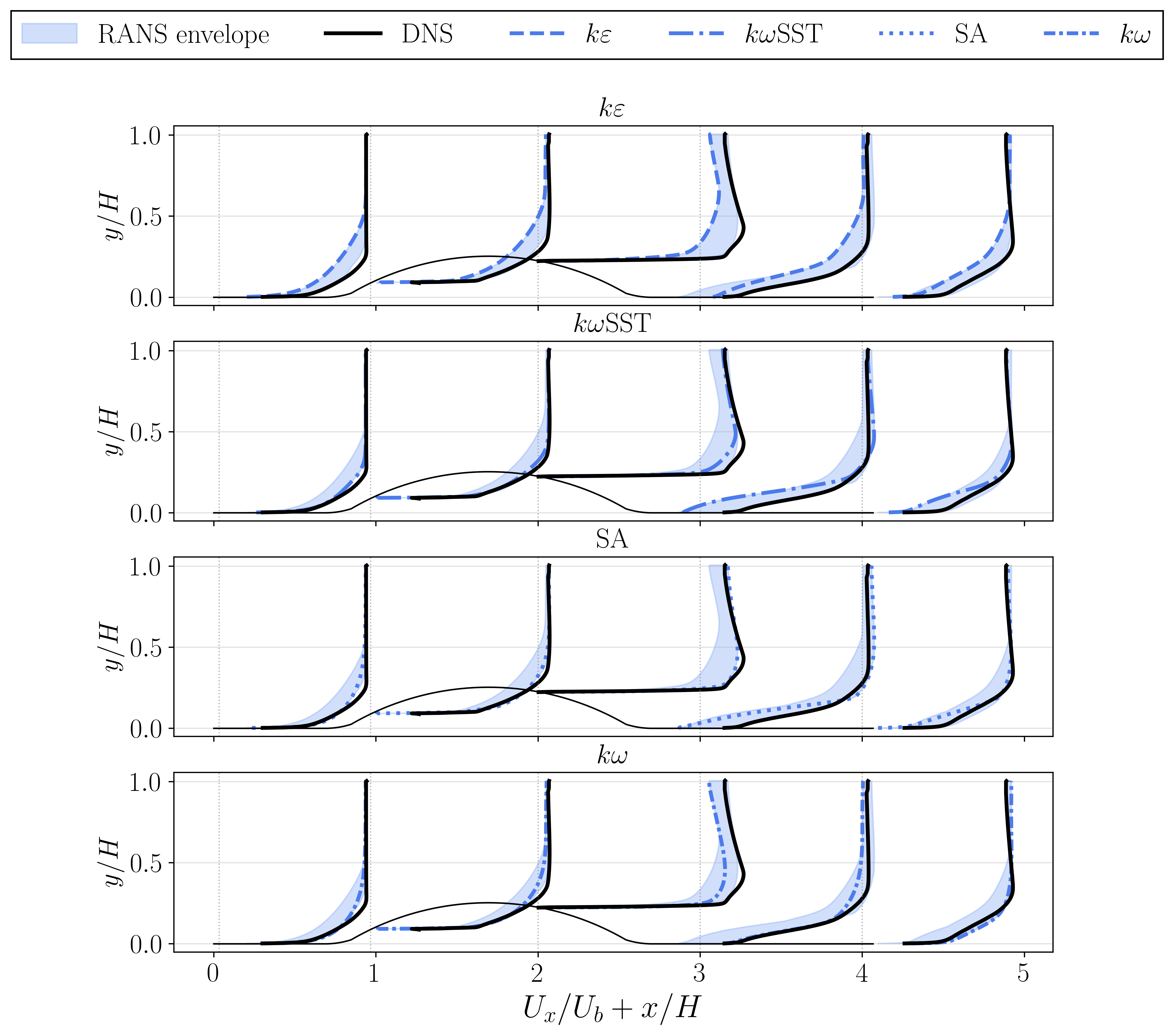}
     \caption{\emph{Test case 2}. Comparison of velocity profiles obtained from RANS simulations using different turbulence models and DNS reference data. We display the horizontal component of the velocity field, extracted at five distinct downstream locations. The results correspond to the test configuration with $h=38\, \mathrm{mm}$. The light blue area indicates the envelope of all the different RANS simulations.}
     \label{Ux_RANS_comp_bumps}
\end{figure}

\par Figure \ref{Ux_RANS_comp_bumps} presents the comparison of the horizontal velocity profiles predicted by the four RANS models against the corresponding high-fidelity data for one representative bump configuration. As in the previous test case, the different turbulence closures exhibit markedly different behaviors. Their predictive accuracy varies substantially across the domain, particularly in regions affected by adverse pressure gradients and mild separation.

As in the previous test case, this pronounced spatial variability further supports the rationale for a locally adaptive model selection strategy. Importantly, the high-fidelity solution is largely contained within the envelope defined by the four RANS predictions, indicating that suitable convex combinations of the baseline models can leverage their complementary features and potentially achieve enhanced accuracy.

\subsubsection{Aggregated models}\label{sec:aggregation-res-bumps}

\par This Subsection presents the results obtained with the aggregated models generated by the two proposed pipelines for the parametrized bump configuration. 
The model performance is assessed on the test set in terms of predictive accuracy, thereby evaluating the capability of the surrogate mixtures to generalize to unseen geometrical configurations and flow conditions. The predictions are compared against the corresponding high-fidelity solutions as well as against the baseline RANS models, allowing for a direct quantification of the improvements introduced by the aggregation framework.

In the KNN aggregation approach, we consider $K=2$ neighbors, as we only have $3$ training solutions.

\begin{remark}
    It is worth highlighting that in the MFR pipeline, after training the aggregation with $N_{\mathrm{train}}^{(mix)}=3$ DNS solutions, we generate $N_{\mathrm{train}}=10$ aggregated solutions for the Mixed-FOM dataset by using the KNN (or the ANN) in inference. Such snapshots are then used to train the ROM. However, in the first test case, we do not need to employ the aggregation strategies in inference in the MFR framework. 
\end{remark}
\begin{figure}[!htpb]
     \centering
     \includegraphics[width=0.75\linewidth]{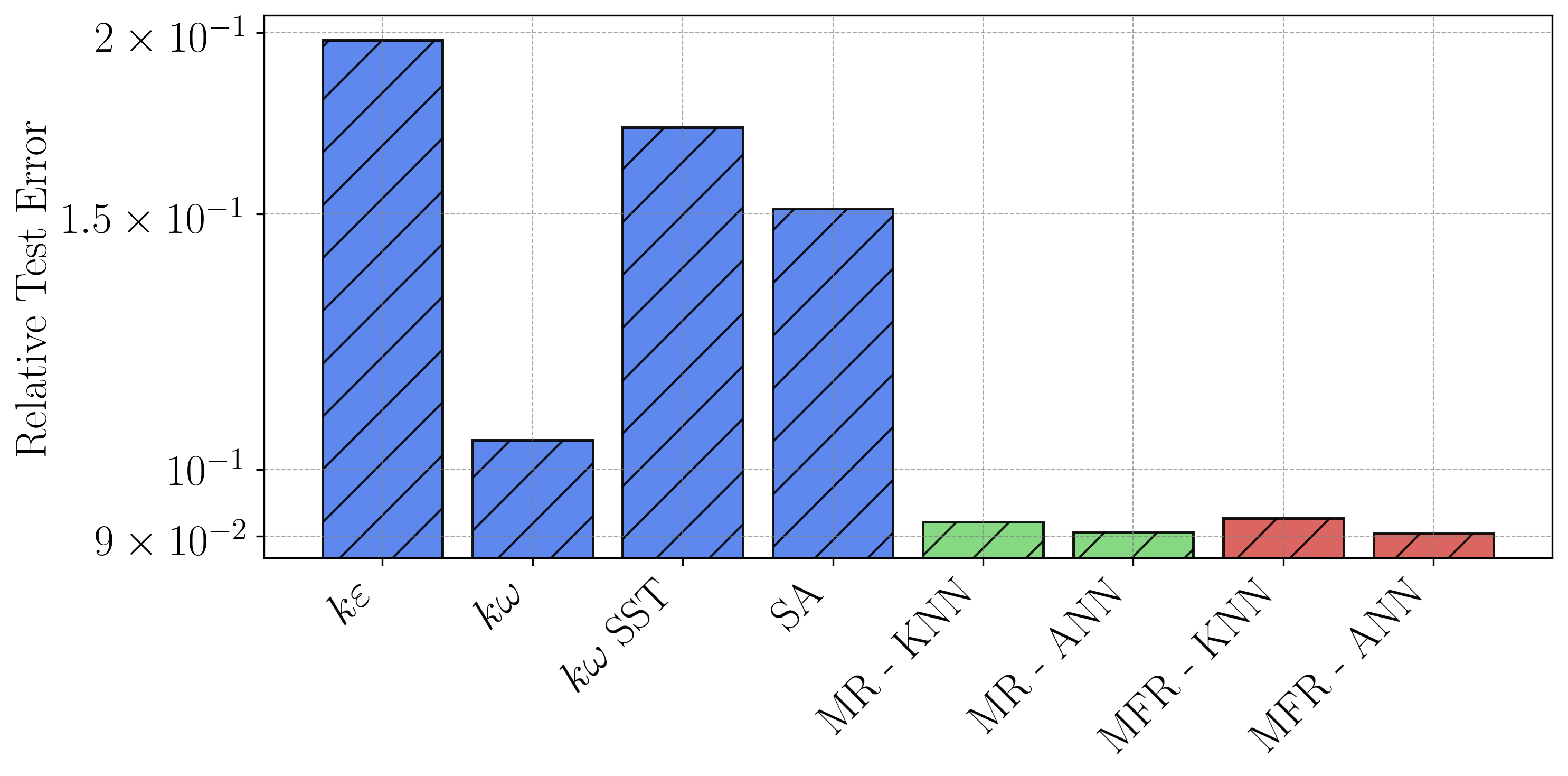}
     \caption{\emph{Test case 2}. Comparison of relative test errors with respect to DNS data for baseline RANS models and surrogate model mixtures.}
        \label{Relative_errors_bumps}
\end{figure}

\par As an initial quantitative assessment, Figure \ref{Relative_errors_bumps} reports the test errors of the individual RANS baseline models together with those of all aggregated surrogate models, computed relative to the high-fidelity reference data. The adopted metric is the mean relative error in the $L^2$ norm. Consistently with the first test case, all aggregated models, regardless of the weighting strategy or the specific aggregation pipeline, outperform the individual RANS component models.

\begin{figure}[!htpb]
     \centering
     \includegraphics[width=0.75\textwidth]{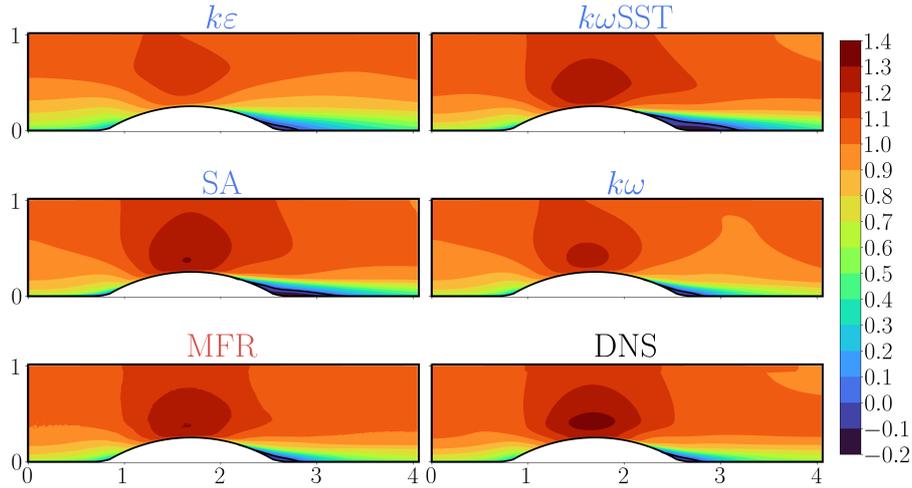}
     \caption{\emph{Test case 2}. Comparison of the horizontal velocity component predicted by the individual RANS models and the MFR-ANN surrogate, with DNS reference solution. The results refer to the \texttt{T2} test configuration.}
        \label{preds_MFROM_bumps}
\end{figure}

\par
\begin{figure}[!htpb]
     \centering
     \includegraphics[width=0.75\textwidth]{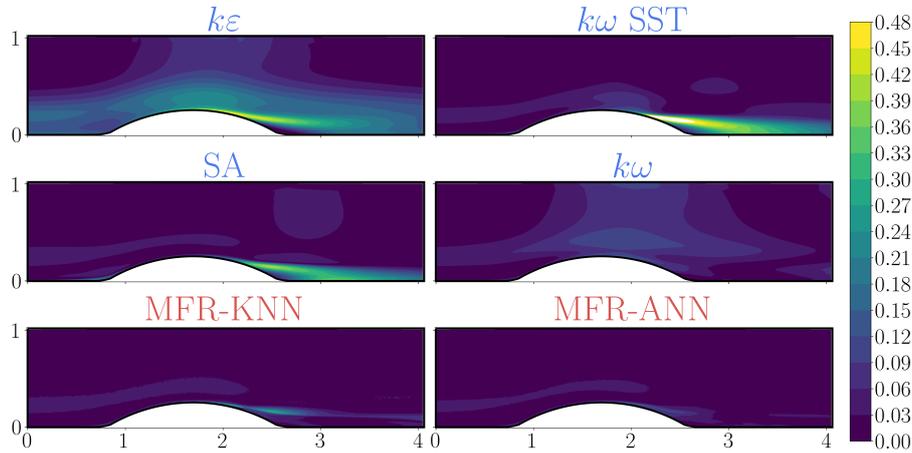}
     \caption{\emph{Test case 2}. Absolute error maps with respect to the DNS benchmark for the horizontal velocity component, predicted by the individual RANS models, the MFR-KNN and the MFR-ANN surrogates. The results correspond to the \texttt{T2} test configuration.}
        \label{ABS_ERRORS_MFROM_bumps}
\end{figure}

 Figure \ref{preds_MFROM_bumps} compares the horizontal velocity fields predicted by the individual RANS models, the MFR-ANN surrogate and the DNS. In particular, the results refer to the \texttt{T2} test configuration, which coincides with $h=38 \, \mathrm{mm}$. This figure highlights the limitations of the individual RANS models. For instance, the recirculation bubble is inaccurately predicted by the 
$k\omega$ SST and SA models. Conversely, the $k\varepsilon$ and $k\omega$ models yield a substantially thicker boundary layer both upstream and downstream of the bump. The MFR prediction, instead, appears overall the closest to the reference solution.
In addition, Figure \ref{ABS_ERRORS_MFROM_bumps} presents the corresponding absolute error maps. These errors are computed with respect to the DNS solution for the individual RANS models as well as the MFR surrogates.
Like in the first test case, the ANN aggregation approach outperforms the KNN, as can be seen both from the global test relative errors (Figure \ref{Relative_errors_bumps}) and from the graphical absolute errors in Figure \ref{ABS_ERRORS_MFROM_bumps}.

\par
\begin{figure}[!htpb]
   \centering
   \begin{subfigure}[b]{0.3\linewidth}
       \centering
       \includegraphics[width=\linewidth]{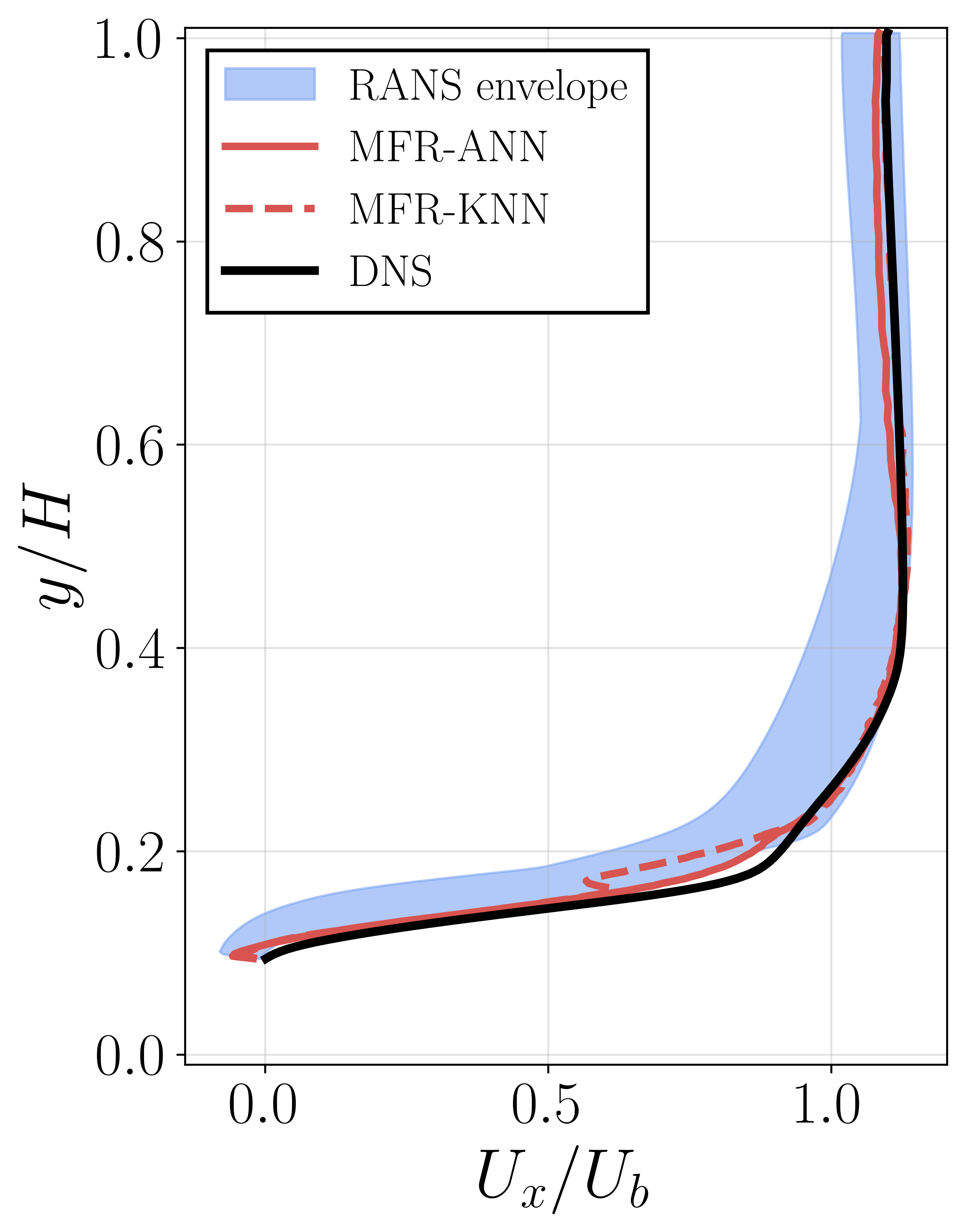}
       \caption{$U_x$ at $x/H=2.4$.}
       \label{profiles_MFROM_U_1_bumps}
   \end{subfigure}
   \hfill
    \begin{subfigure}[b]{0.3\linewidth}
        \centering        
        \includegraphics[width=\linewidth]{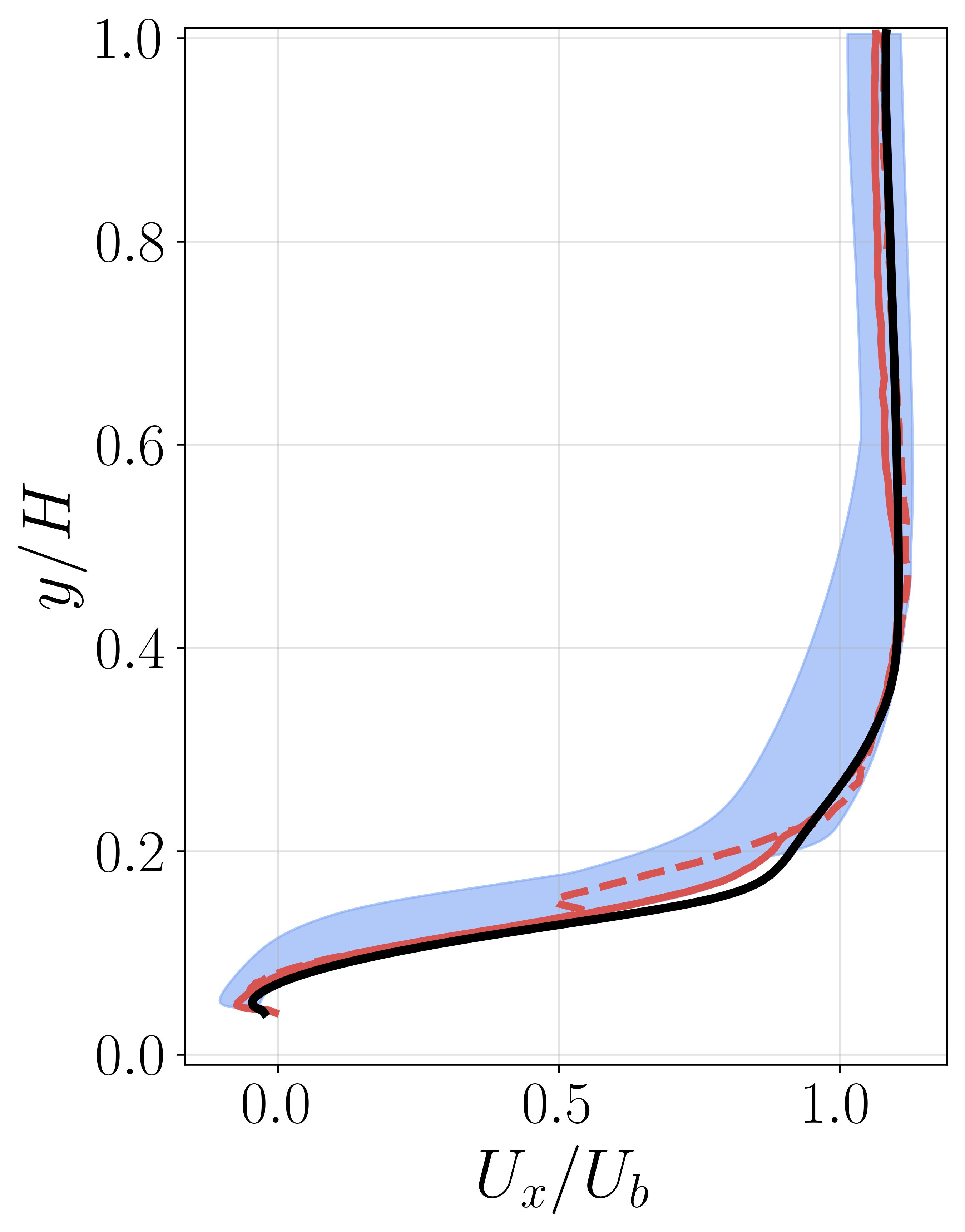}
        \caption{$U_x$ at $x/H=2.5$.}
        \label{profiles_MFROM_U_2_bumps}
    \end{subfigure}
   \hfill
    \begin{subfigure}[b]{0.3\linewidth}
        \centering        
        \includegraphics[width=\linewidth]{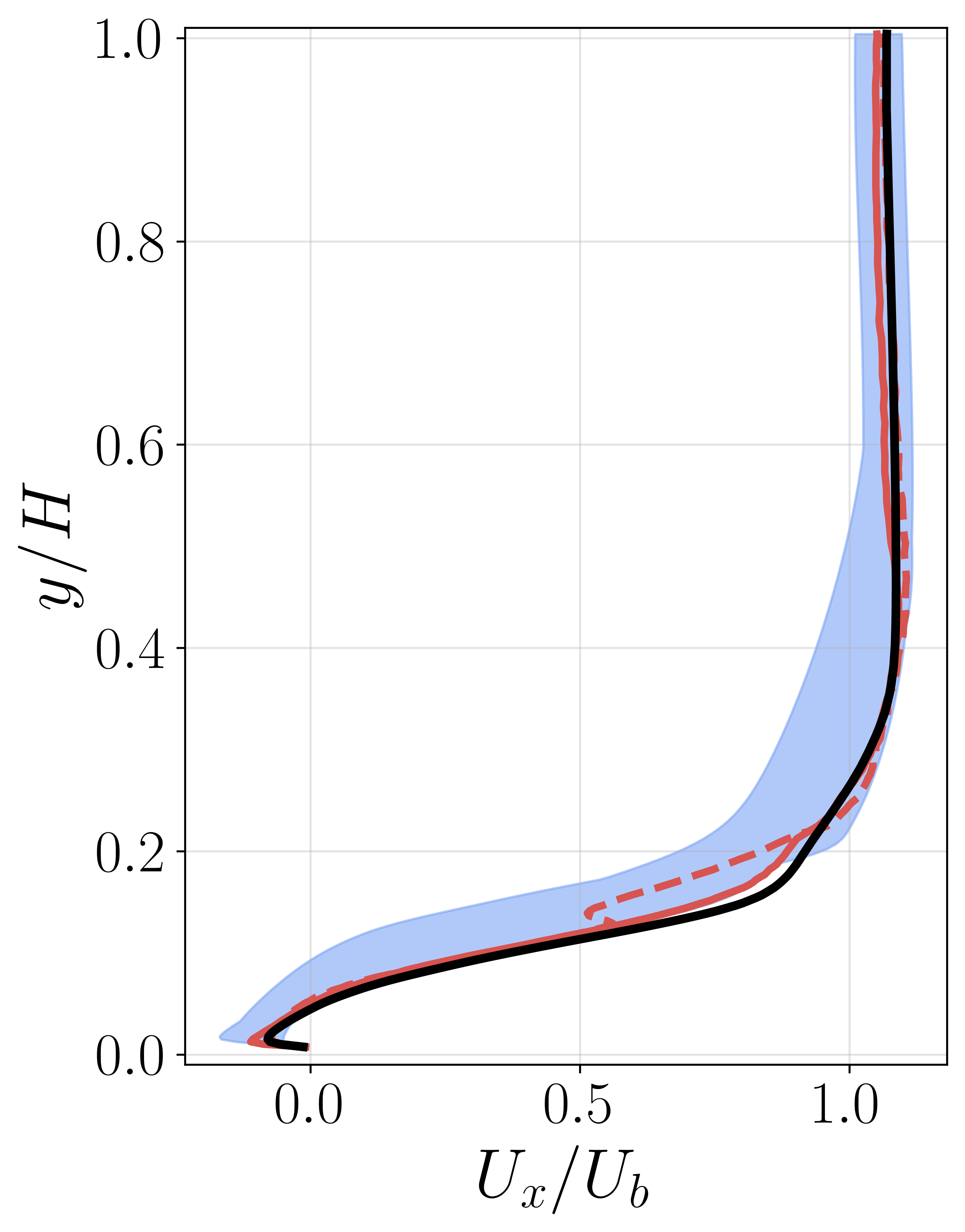}
        \caption{$U_x$ at $x/H=2.6$.}
        \label{profiles_MFROM_U_3_bumps}
    \end{subfigure}
    \caption{\emph{Test case 2}. Comparison of horizontal velocity profiles obtained from the two MFR models employing different weighting techniques, alongside DNS reference data. The profiles are extracted in the separation region at $x/H=2.4$ (left), $x/H=2.5$ (center) and $x/H=2.6$ (right) for the \texttt{T2} test configuration. The shaded region indicates the prediction range of the four individual RANS turbulence models.}
    \label{profiles_MFROM_U_bumps}
\end{figure}

 Figure \ref{profiles_MFROM_U_bumps} shows that both surrogate models produce smooth solutions that are fairly consistent with each other. The predicted velocity profiles generally agree with the DNS data, capturing the overall reference behavior despite the considerable variability among the individual RANS component models, as indicated by the wide shaded region. Notably, the ANN-based approach now provides only a mild improvement over the KNN-based model, which still tends to show the largest deviations, in line with the trends observed in the previous analysis.

\par
\begin{figure}[!htpb]
    \centering
    \begin{subfigure}[b]{\linewidth}
        \centering        
        \includegraphics[width=0.75\textwidth]{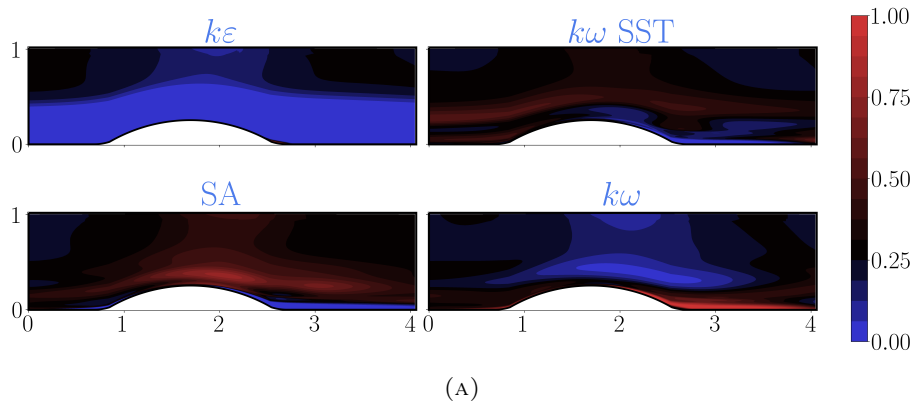}
        \caption{}
        \label{weightsKNN_U_bumps}
    \end{subfigure}
    \caption{\emph{Test case 2}. Spatial distribution of the weighting functions for the four component RANS models, in the MFR-KNN pipeline.}
    \label{weightsKNN_MF_bumps}
\end{figure}

 Figure \ref{weightsKNN_MF_bumps} shows the spatial distributions of the weighting functions assigned to the four component RANS turbulence models, as determined by the MFR-KNN approach.

\par The space-dependent weights from the KNN strategy indicate very low contributions from the $k\varepsilon$ model within the boundary layer, strong activation of the $k\omega$ model near the recirculation bubble, and dominant contributions from the SA and $k\omega$ SST models near the top of the bump. Overall, the weight distributions continue to reflect the local predictive capabilities of each turbulence closure, demonstrating that the KNN approach effectively identifies the best-performing RANS model in each region of the domain.

\par
\begin{figure}[!htpb]
     \centering
     \includegraphics[width=0.75\textwidth]{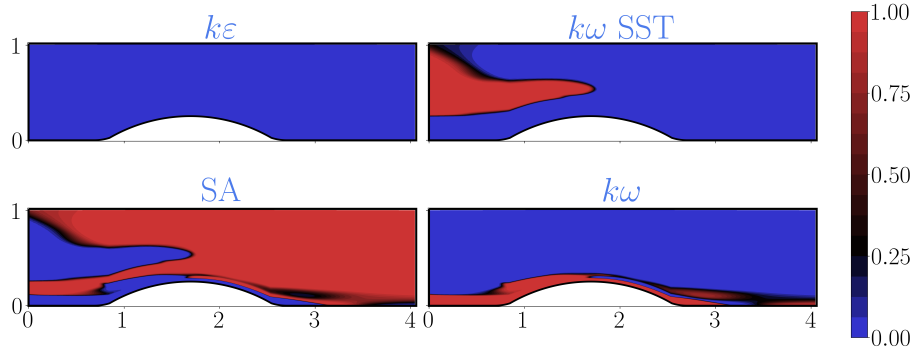}
     \caption{\emph{Test case 2}. Spatial distributions of the weighting functions for the four component RANS models, in the MFR-ANN pipeline.}
        \label{weightsANN_MF_bumps}
\end{figure}

 Similarly, Figure \ref{weightsANN_MF_bumps} illustrates the spatially varying weights predicted by the ANN-based weighting approach. As observed in the first case, the ANN weight distributions show substantially less variability compared to those from the KNN method. The ANN approach tends to assign the majority of the weight to a single best-performing model over large regions of the domain, with the other closures contributing minimally or being entirely neglected. As a result, the representation of model uncertainty is limited, and the weighting patterns are less interpretable. Nonetheless, similar trends to the KNN strategy are evident: the $k\varepsilon$ model is largely disregarded throughout, the $k\omega$ model dominates in the recirculation zone, and the SA model primarily governs the region near the top of the bump.
\begin{remark} Also in this case, we remark that the space-dependent weights are used only to build the high-fidelity aggregated predictions that constitute the Mixed-FOM training dataset. During the online phase, the surrogate models do not actually employ these weight fields, which are provided exclusively to show the effectiveness of the weighting strategies.
\end{remark}

\subsection{Computational costs}\label{sec:computational-costs}

\par This Subsection examines the computational costs associated with both the offline and online phases of the proposed aggregation pipelines.
We report here the times regarding the first test case, but similar observations can be drawn for the second test case.

In particular, we show that, leveraging the ROMs, the surrogates we developed achieve very efficient solutions in their online phase and thus provide significant speed-ups compared to the standard RANS simulations.

\subsubsection{Costs of the offline stage}
\par The offline phase is performed only once for each surrogate mixture and entails the most computationally demanding steps of the proposed aggregation pipelines. Specifically, this stage is composed of the following main tasks:
\begin{enumerate}
    \item \textbf{Execution of the baseline RANS simulations}. A total of $N_{\mathrm{train}}^{\mathrm{tot}}$ baseline RANS simulations are performed, serving as the foundation of the proposed methodologies. In particular, the dataset of high-fidelity snapshots is generated only once and is subsequently reused for the construction of all surrogate mixtures developed within both pipelines. Each FOM simulation requires, on average, approximately $1200 \mathrm{s}$ of CPU time when run on a single processor core of the SISSA HPC cluster \textit{Ulysses} (200 TFLOPS, 2TB RAM, 7000 cores).
    \item \textbf{Building the non-intrusive reduced order models}. In the first aggregation pipeline a single ROM is employed, whereas in the second pipeline four distinct ROMs are constructed. In both cases, the ROMs use an autoencoder for dimensionality reduction combined with RBF interpolation for the approximation stage. The computational cost of this step is dominated by the AE training time, which is approximately $90 \mathrm{s}$ on an Apple M2 Pro CPU (16 GB RAM), utilizing all available $12$ processor cores in parallel. It is worth noting that, even in this offline phase, the time required to train an AE is substantially lower than that of a single RANS simulation, which indicates that in this context the additional overhead introduced by the ROM framework is negligible compared to the cost of generating the high-fidelity data.
    \item \textbf{Training of the space-dependent weighting technique}. The computational cost of this step varies significantly depending on the weighting strategy considered:
    \begin{itemize}
        \item \textbf{KNN supervised regression}: the total cost is approximately $3 \mathrm{s}$ on an Apple M2 Pro CPU (16 GB RAM) using a single processor core. The majority of this time is spent on solving the minimization problem required to estimate the $\varsigma$ hyperparameter.
        \item \textbf{ANN approach}: a feed-forward artificial neural network is trained, requiring $\sim1800 \mathrm{s}$ when the computations are offloaded to an Apple M2 Pro 19-cores GPU.
    \end{itemize}
\end{enumerate}

\subsubsection{Costs of the online stage}
\par During the online phase, the trained surrogate models are employed to generate predictions for new values of the problem parameters. Since this process is typically executed a very large number of times, it is essential for this stage to remain computationally cheap and substantially faster than the baseline full-order model simulations. All numerical tests performed for this part of the analysis are carried out on a single core of an Apple M2 Pro CPU (16 GB RAM). The corresponding results are reported in Table \ref{speedups}.
\par In the case of the Mixed FOM--ROM aggregation pipeline, the online stage consists solely of a single ROM evaluation, which requires approximately $4\times10^{-4}\ \mathrm{s}$. This corresponds to a remarkable speed-up of the order of $10^{6}$ relative to the average CPU time needed for a baseline RANS simulation.
\vspace{0.1cm}
\begin{table}[!htpb]
    \caption{\emph{Test case 1}. Online execution times of the surrogate models and corresponding speed-up factors with respect to the baseline RANS simulations.}
    \centering
    \begin{tabular}{c c c}
    \toprule
         \makecell[c]{\textbf{Aggregation}\\\textbf{pipeline}} & \makecell[c]{\textbf{Online time $[\mathrm{s}]$}} & \makecell[c]{\textbf{Speed-up}\\\textbf{factor}}\\
         \midrule
         MFR  & $\sim4\times10^{-4}$ & $\mathcal{O}(10^{6})$\\
         MR  & $\sim6\times10^{-3}$ & $\mathcal{O}(10^{5})$\\
    \bottomrule
    \end{tabular}
    \label{speedups}
\end{table}
\par For the Mixed-ROM pipeline, the online phase is slightly more complex. Specifically, it includes the computation of four ROM solutions, one for each component model, followed by the evaluation of the space-dependent weighting strategy and the assembly of the final aggregated prediction. Despite the additional steps, the procedure remains highly efficient, with a total runtime of approximately $6\times10^{-3}\ \mathrm{s}$, thus still yielding a significant speed-up of the order of $10^{5}$ compared to the baseline FOM simulations. Importantly, the choice of the weighting technique does not produce any noticeable difference in the online execution times.

\section{Conclusions}
\label{sec:conclusions}
In this work, we have proposed a unified and data-driven framework that combines space-dependent aggregation of multiple RANS turbulence models with non-intrusive nonlinear reduced order modeling. The objective was to simultaneously enhance predictive accuracy and achieve computational efficiency suitable for many-query and near real-time applications.

The first novelty is the introduction and comparison of two distinct aggregation strategies: the Mixed FOM–ROM (MFR) pipeline, where a single ROM is trained on aggregated high-fidelity RANS solutions, and the Mixed-ROM (MR) pipeline, where separate ROMs are first constructed and then aggregated at the reduced level. This comparison clarifies the trade-offs between offline complexity and overall efficiency. Although the MR approach is conceptually appealing, as it aggregates fully independent reduced models, it does not provide a significant accuracy advantage over MFR, while being more expensive in the offline stage. For this reason, the MFR pipeline emerges as the most effective compromise between simplicity, cost, and predictive performance.

The second important novelty lies in the learning of the space-dependent aggregation weights. In addition to a classical KNN-based regression strategy, we introduced an ANN-based formulation that directly outputs spatially continuous and normalized weights. The ANN approach consistently demonstrated improved predictive performance on unseen configurations in both test cases, while removing the need for hyperparameter tuning associated with kernel-based methods.

The methodology was validated on two parametrized benchmark problems characterized by separated turbulent flows. In both cases, the aggregated surrogate models significantly outperformed the individual RANS closures and their corresponding ROM counterparts. At the same time, the proposed surrogates achieved speed-ups of up to six orders of magnitude with respect to baseline RANS simulations. These results confirm that the framework successfully combines improved physical accuracy with drastic computational acceleration.

Overall, this study demonstrates that the integration of space-dependent aggregation and nonlinear non-intrusive ROMs provides an efficient and robust pathway toward accurate surrogate modeling of turbulent flows.

Several directions naturally follow from this work. Extending the framework to three-dimensional configurations and higher Reynolds numbers would further assess its scalability in industrially relevant scenarios. The integration within multi-fidelity environments, combining sparse DNS or LES data with large RANS datasets, represents another promising avenue. Finally, incorporating uncertainty-aware weighting strategies could improve reliability assessment in predictive settings.

\section*{Acknowledgments}
This study was funded by the European Union - NextGenerationEU, in the framework of the iNEST - Interconnected Nord-Est Innovation Ecosystem (iNEST ECS00000043 – CUP G93C22000610007). The views and opinions expressed are solely those of the authors and do not necessarily reflect those of the European Union, nor can the European Union be held responsible for them. The authors would like to sincerely thank Paola Cinnella for her valuable input and insightful discussions, which greatly contributed to improving the manuscript.

\appendix
\section{Supplementary material}\label{sec:sup}
\par In this supplementary Section, we provide the description of the architectures and the training hyperparameters of the neural networks employed throughout this study, to ensure reproducibility. In particular, Table \ref{ANN_details} reports the specifications of the feed-forward fully connected ANNs used to compute the optimal space-dependent weights of the model mixtures adopting the ANN weighting technique. The same configuration is selected for all such ANNs for consistency. Conversely, Table \ref{AE_table} presents the details of the autoencoders employed in the dimensionality reduction stages of all ROMs constructed in this work, for which a common set of settings is used. Notably, for each class of networks, the same specifications were considered across all instances in order to ensure a fair comparison and to avoid potential biases in the assessment of the results.
\begin{table}[!htpb]
    \caption{Architecture and training parameters of the artificial neural networks employed within the ANN space-dependent weighting strategy.}
    \centering
    \begin{tabular}{c c c c c c}
    \toprule
         \textbf{Test case}&\makecell[c]{\textbf{Hidden}\\\textbf{layers}} & \makecell[c]{\textbf{Activation}\\\textbf{function}} & \makecell[c]{\textbf{Learning}\\\textbf{rate}} & \makecell[c]{\textbf{Weight}\\\textbf{decay}} & \textbf{Epochs} \\
         \midrule
         1&$[30,30,30,30,30]$ & Softplus & $5\times10^{-4}$ & $1\times10^{-4}$ & $60000$\\
         2&$[50, 50, 50]$ & Tanh & $1\times10^{-3}$ & $1\times10^{-5}$ & $10000$\\
    \bottomrule
    \end{tabular}
    \label{ANN_details}
\end{table}

\begin{table}[!htpb]
    \caption{Autoencoders architecture and training parameters.}
    \centering
    \begin{tabular}{c c c c c c}
    \toprule
         \textbf{Test case}&\makecell[c]{\textbf{Hidden}\\\textbf{layers}} & \makecell[c]{\textbf{Activation}\\\textbf{function}} & \makecell[c]{\textbf{Learning}\\\textbf{rate}} & \makecell[c]{\textbf{Weight}\\\textbf{decay}} & \textbf{Epochs} \\
         \midrule
   1& $[50,20,3,3,20,50]$ & Softplus & $5\times10^{-4}$ & $1\times10^{-4}$ & $10000$\\
   2& $[50,20,10,10,20,50]$ & Softplus & $1\times10^{-3}$ & $1\times10^{-8}$ & $10000$\\
    \bottomrule
    \end{tabular}
    \label{AE_table}
\end{table}

\bibliographystyle{abbrv}
\bibliography{biblio}

\end{document}